\documentclass[draftcls,onecolumn,12pt]{IEEEtran}

\ifCLASSINFOpdf
\else
\fi

\usepackage{setspace}

\usepackage{color}
\usepackage{amsmath,color,amsthm}
\usepackage{amssymb}  
\usepackage{amsfonts}
\usepackage{graphicx}
\usepackage{dsfont}
\usepackage{psfrag}
\usepackage{subfigure,epstopdf}
\usepackage[norefs,nocites]{refcheck}           
\usepackage{cite}

\usepackage[ruled]{algorithm2e}
\usepackage{amsmath}

\newtheorem{theorem}{Theorem}
\newtheorem{corollary}{Corollary}
\newtheorem{remark}{Remark}
\newtheorem{proposition}{Proposition}
\newtheorem{lemma}{Lemma}
\newtheorem{definition}{Definition}

\begin{document}
%

\title{
A  Framework for Structural  Input/Output\\
 and Control Configuration Selection  \\
in Large-Scale Systems}

\author{\IEEEauthorblockN{S\'ergio Pequito $^{\dagger,\ddagger}$ \qquad Soummya Kar $^{\dagger}$ \qquad A. Pedro Aguiar $^{\ddagger,\diamond}$}
\thanks{ $^{\dagger}$ Department of Electrical and Computer Engineering, Carnegie Mellon University, Pittsburgh, PA 15213

$^{\ddagger}$  Institute for System and Robotics, Instituto Superior T\'ecnico, Technical University of Lisbon, Lisbon, Portugal

$^{\diamond}$ Department of Electrical and Computer Engineering, Faculty of Engineering, University of Porto, Porto, Portugal
}}


\setstretch{1.42}

\maketitle

\vspace{-1.7cm}

{\setstretch{1.44}
\begin{abstract}


This paper addresses problems on the structural design of control systems taking explicitly into consideration the possible application to large-scale systems. We provide an efficient and unified  framework to solve the following major minimization problems: (i) selection of the minimum number of manipulated/measured variables to achieve structural controllability/observability of the system, and (ii) selection of the minimum number of  feedback interconnections between measured and manipulated variables  such that the closed-loop system has no structurally fixed modes. Contrary to what would be expected, we show that it is possible to obtain a global solution for each of the  aforementioned minimization problems using polynomial complexity algorithms in the number of the state variables of the system. In addition, we provide several new graph-theoretic characterizations of structural systems concepts, which, in turn, enable us to characterize all possible solutions to the above problems.

\end{abstract}

}

\IEEEpeerreviewmaketitle

\vspace{-0.7cm}

\section{Introduction}\label{secIntro}
{\setstretch{1.45}

This paper is motivated by the dearth of scalable techniques for the analysis and synthesis of  large-scale complex systems, notably ones which tackle design and decision making in a single framework. Examples include power systems, public or business organizations, large manufacturing systems, wireless control systems, biological complex networks, and formation control, to name a few.  A central and challenging issue that arises when dealing with such complex systems is that of structural design. In other words, given a \emph{plant} of a system, we are  interested in providing a framework that addresses the following questions \cite{Skogestad04a}:

\begin{enumerate}
\item[1)] Which variables should be measured?
\item[2)] Which variables should be manipulated?
\item[3)] Which feedback links should be incorporated between the  sets defined in 1) and 2)?
\end{enumerate}
}

Problems 1)-2) are commonly referred to as the  \emph{input/output (I/O) selection problem}, whereas  problem 3) is referred to as the \emph{control configuration (CC) selection problem} \cite{reviewIO}. The latter problem is of significant importance  in the area of decentralized control, where the goal is to understand which subset of sensors (outputs) and local controllers (inputs) need to be feedback connected so that specific properties (e.g., stability) of the overall system hold.
The choice of inputs and outputs affects the performance, complexity and costs of the control system. Due to the combinatorial nature of the selection problem, efficient and systematic methods are required to complement the designer intuition, experience and physical insight \cite{reviewIO}.

 Motivated by the above problems, in this paper we provide an efficient framework, that addresses both  the I/O and CC problems for (possibly large-scale) linear time invariant systems, by resorting to structural systems theory \cite{dionSurvey}, where the main idea is to investigate system-theoretic properties based only on the sparsity pattern (i.e., location of zeroes and non-zeroes) of the system matrices. Structural systems based formulations offer the added advantage of being able to deal with scenarios in which the specific numerical values of the system parameters are not accurately known.   The major design constraints or system properties that  are addressed in the  I/O selection problem are those of controllability and observability, which in the context of structural systems correspond to  \emph{structural controllability}  and \emph{structural observability} (to be formally defined in Section II). In addition, the absence of  \emph{structurally fixed modes} is the key property of interest   to ensure in the CC selection problem due to its implications on  \emph{generic} pole placement for decentralized control systems \cite{sezerFixedModes}. Design and analysis based on structural systems provide  system-theoretic guarantees that hold for almost all numerical instances of the parameters, except on a manifold of zero Lebesgue measure \cite{Reinschke:1988}.

 We now describe precisely the problems addressed in this paper.

\vspace{-0.55cm}
\subsection*{Problem Statement}

\vspace{-0.2cm}

 Consider a given (possibly large-scale) plant with autonomous dynamics
\vspace{-0.4cm}
\begin{equation}
\dot x = A x,
\vspace{-0.4cm}
\label{dynSys}
\end{equation}
where $x \in\mathbb{R}^n$ denotes the state of the plant and $A$ is an $n\times n$ matrix. Suppose that the sparsity (i.e., location of zeroes and non-zeroes) pattern of $A$ is available, but the specific numerical values of its non-zero elements are not known. Let $\bar A \in \{0,1\}^{n\times n}$ be the binary matrix that represents the structural pattern of $A$, i.e., it encodes the sparsity pattern of $A$  by assigning $1$ to each non-zero entry of $A$ and zero otherwise.

$\rhd$ \underline{Sparsest I/O selection problem}

$\mathcal{P}_1$ Given $\bar A$ associated with \eqref{dynSys}, find structural input and output matrices $(\bar B,\bar C)$ that solve
\vspace{-.3cm}
\begin{equation}
\underbrace{
\begin{array}{cc}
\min\limits_{\bar B\in \{0,1\}^{n \times n}}  &  \|\bar B\|_0\\
\text{s.t.} & (\bar A, \bar B) \text{ is struct. controllable}
\end{array}}_{(2a)}
\begin{array}{c}
+\\
\\
\end{array}
\underbrace{\begin{array}{cc}
\min\limits_{\bar C\in \{0,1\}^{n \times n}}  &  \|\bar C\|_0\\
\text{s.t.} & (\bar A, \bar C) \text{ is struct. observable}
\end{array}}_{(2b)}
\label{p1}
\end{equation}
where $\|M\|_0$ is the semi-norm that denotes the number of non-zero entries in the binary matrix $M\in \{0,1\}^{n\times n}$ and $(\bar A,\bar B)$ (respectively, $(\bar A,\bar C)$) is the pair of matrices that represent the structural system with dynamics matrix structure $\bar A$ and input (respectively output) matrix structure $\bar B$ (respectively $\bar C$).
\hfill $\diamond$

Problems (2a) and (2b) correspond to the sparsest input and output selection problem, respectively.  
Note that a solution to $\mathcal P_1$ may not necessarily be one that will correspond to an implementation with a \emph{minimum number of inputs/outputs}. In the paper, we will also characterize the subset of all solutions to $\mathcal P_1$ which have that property. In other words, we determine the sparsest input/output matrix, comprising the smallest number of inputs/outputs, and ensuring structural controllability/observability. Additionally, we are interested in obtaining solutions to more constrained variants of $\mathcal{P}_{1}$, specifically, that of characterizing structurally controllable/observable configurations with the \emph{minimum number of dedicated inputs/outputs}. Dedicated input configurations are those in which each input may only manipulate a single state variable (i.e., at most one entry in each column of $\bar B$ can be non-zero), whereas,  dedicated output configurations are those in which each output corresponds to a sensor measuring a single state variable (i.e., at most one entry in each row of $\bar C$ can be non-zero). Formally, this last problem can be posed as follows:

$\mathcal{P}_1^d$  Given $\bar A$ associated with \eqref{dynSys},  find $(\bar B,\bar C)$ that solve
\vspace{-.3cm}
\begin{equation}
\underbrace{
\begin{array}{cc}
\min\limits_{\bar B\in \{0,1\}^{n \times n}}  &  \|\bar B\|_0\\
\text{s.t.} & (\bar A, \bar B) \text{ is struct. controllable}\\
& \|\bar B_{.j}\|_0\le 1, \ j=1,\cdots,n
\end{array}}_{(3a)}
\begin{array}{c}
+\\
\\
\end{array}
\underbrace{\begin{array}{cc}
\min\limits_{\bar C\in \{0,1\}^{n \times n}}  &  \|\bar C\|_0\\
\text{s.t.} & (\bar A, \bar C) \text{ is struct. observable}\\
& \|\bar C_{i.}\|_0\le 1, \ i=1,\cdots,n
\end{array}}_{(3b)}
\label{p1d}
\vspace{-0.25cm}
\end{equation}
where $\bar B_{.j}$ represents the $j$th column of $\bar B$ and $\bar C_{i.}$ the $i$th row of $\bar C$.
\hfill $\diamond$

Interestingly, by formulation, $\mathcal{P}_{1}^{d}$ (initially addressed in [6]) appears to be more constrained than $\mathcal{P}_{1}$, however, we will show that solutions to $\mathcal{P}_{1}^{d}$ constitute a subclass (generally strict) of solutions to $\mathcal{P}_{1}$, upon which the general solution to $\mathcal P_1$ can be described.

To formally state the  CC selection problem, for a linear system $(A,B,C)$ with inputs and outputs, let  $\bar K$  represent the \emph{information pattern}, i.e., $\bar K_{ij}=1$ if output $j$ is available to input $i$ and $0$ otherwise. The existence of feedback matrices $K$ with the same sparseness of $\bar K$ that allow arbitrary pole-placement of the  output feedback closed-loop system (given by $A+BKC$), is associated with the notion of \emph{fixed modes} of the system $(A,B,C)$  w.r.t. (with respect to) the  specified information pattern $\bar K$ \cite{Davison}. Its structural counterpart is the notion of  \emph{structurally fixed modes} \cite{sezerFixedModes}: given $(\bar A,\bar B,\bar K,\bar C)$, we say that $(\bar{A},\bar{B},\bar{C})$ has no structurally fixed modes w.r.t. $\bar{K}$, if for almost all realizations of $(A,B,C)$ (with the same structure as $(\bar A,\bar B,\bar C)$) there exist feedback matrices $K$ (with the same structure as $\bar K$) such that  the poles of the corresponding   static output feedback closed-loop system  can be  arbitrarily placed in the complex plane (and in particular in the left half  of the complex plane). In addition, we further notice that  the use of dynamic compensators  is accounted for in the scope of static output feedback to achieve a specified pole-placement. However, in general, the  use of dynamic gain matrices fall beyond the scope of the current study, which would lead to additional design flexibility but might be prohibitive to compute in the setting of large-scale dynamical systems. Consequently, the absence of structurally fixed modes is the  property of interest  that we seek to ensure in the (CC) selection problem. Often, deciding which outputs are to be used  by each of the inputs is called the partitioning or pairing problem, when the information pattern is additionally restricted to have a  block diagonal structure \cite{reviewIO}. In this paper, we focus on a  more general problem where the information pattern is not restricted to be  block diagonal. 

Formally, the maximum sparseness jointly  I/O and CC selection problem is stated as follows:
\vspace{-0.2cm}

$\rhd$ \underline{Jointly sparsest I/O and CC selection problem}

\begin{itemize}
\item[$\mathcal{P}_2$]  Given $\bar A$ associated with \eqref{dynSys},   find $(\bar B,\bar K,\bar C)$ that solve
\vspace{-.4cm}
\begin{align}
\min_{\bar B,\bar K,\bar C\in \{0,1\}^{n \times n}} &\quad  \|\bar B\|_0+ \|\bar C\|_0 +\|\bar K\|_0 \label{p2}\\
\text{s.t.} \qquad & \quad (\bar A,\bar B,\bar K, \bar C) \text{ has no struct. fixed modes}
\notag
\end{align}
\end{itemize}
\vspace{-1cm}
\hfill $\diamond$

The usual approach to solve the I/O and CC selection problem (often to achieve different goals, other than maximum sparseness) is to address them  independently and sequentially \cite{reviewIO}, but with no guarantee that such solution yields the optimal. In the present paper, we actually show that \emph{some} solutions of the I/O selection problem $\mathcal P_1$ can be used to solve the joint I/O and CC selection problem $\mathcal P_2$ in an optimal fashion. Conversely, we demonstrate that all possible solutions to $\mathcal{P}_{2}$ may be characterized in terms of solutions $(\bar{A},\bar{B},\bar{C})$ of $\mathcal{P}_{1}$ together with a construction referred to as \emph{mixed-pairing} (informally a pairing between inputs and outputs).

\vspace{-0.4cm}

{\setstretch{1.45}
\subsection*{Related Work}
\vspace{-0.2cm}

Both the I/O and CC selection problems have received significant attention  in the literature, see \cite{reviewPlacement,reviewIO} and references therein. In the context of the current work, we restrict our attention to those papers which study the above problems in the structural systems framework. Since the seminal paper  \cite{Lin_1974} in structural systems theory, a large number of papers have considered several variants of the I/O and CC selection problems in which different solution criteria and applications are presented,  see \cite{SiljakBook,largeScale,Reinschke:1988,Murota:2009:MMS:1822520} and references therein (see also  \cite{dionSurvey} for a very useful survey of several important results in structural systems theory). For instance, in   \cite{CommaultD13} and references therein, given  the dynamics structure of a linear-time invariant system dynamics and a possible collection of inputs, the objective was to determine the minimum subset of inputs that achieve structural controllability, which we refer to as the \emph{constrained minimal input selection problem}.  Similar work is presented in \cite{DBLP:journals/automatica/BoukhobzaH11a}, where the analysis of I/O selection is foreseen, yet in a more general setting: that of ensuring structural observability for linear systems in descriptor form with unknown inputs. Determining feasible solutions to the above constrained minimal input selection problem has beena major focus of the structural systems literature, see  \cite{dionSurvey,Murota:2009:MMS:1822520,DBLP:journals/automatica/BoukhobzaH11a} for representative work. A large majority of the proposed solution methods to the above constrained minimal input selection problem rely on a two-step optimization procedure (see~ \cite{CommaultD13}) which, in general, leads to suboptimal solutions as recently emphasized in~ \cite{CommaultD13}. On a related note, the constrained minimal input selection problem was in fact shown to be NP-complete in general, and, hence are unlikely to have polynomial complexity algorithmic solutions. In contrast the I/O selection variant that we study in this paper is quite different: more precisely, we aim to determine the sparsest I/O matrices that ensure structural controllability/observability.  Moreover, the sparsest I/O selection problem is shown to be polynomially solvable in this paper, which is achieved through design techniques that are very different from the ones developed in~\cite{dionSurvey,Murota:2009:MMS:1822520,DBLP:journals/automatica/BoukhobzaH11a} for addressing the constrained minimal I/O selection problem.

The problem of identifying the minimum number of inputs/outputs required to achieve structural controllability/observability was considered in \cite{dionSurveyKyb,liu11}. Specifically, in the context of complex networks, it was shown in~\cite{liu11} that the minimum number of controlling agents required to achieve structural controllability is related to the number of right unmatched vertices of an associated bipartite matching problem. However, in general (specifically, when the system digraph consists of multiple strongly connected components), such characterization of the minimum number of inputs/outputs required to achieve structural controllability/observability is not sufficient to address the sparsest I/O selection problems $\mathcal{P}_{1}^{d}$ and $\mathcal{P}_{1}$ considered in this paper, since the latter problems additionally require identifying the minimum number of connections that are to be made between the state variables and the control inputs (see also remarks after Corollary 1 for a more technical discussion). This latter characterization is obtained in this paper which enables us to address $\mathcal{P}_{1}$ and $\mathcal{P}_{1}^{d}$ in full generality. In this context, we also note \cite{Pequito}, our preliminary work, which provides the first general solution to the dedicated input/output assignment problem $\mathcal{P}_{1}^{d}$ and further provides polynomial complexity algorithms to explicitly compute such solutions. On a related note, more recently, in \cite{mincontrolNP} it was shown that  the minimal controllability problem (sparsest input design to ensure controllability given a numerical instance, instead of structural controllability) is NP-hard (w.r.t. the system size). Interestingly, in this paper, we show that the structural counterpart of this problem, i.e., problem (2a) (see (2)), is of polynomial complexity in the system size.

 In \cite{Khan,Khan2}, the design of a network of sensors and communication among those are sought  to ensure sufficient conditions for distributed state estimation with bounded error. Design of networked control systems is pursued in \cite{Sundaram,DBLP:journals/jsac/PajicMPS13}, where given a decentralized plant, modeled as a discrete linear time invariant system equipped with actuators and sensors, the communication topology design between actuators and sensors to achieve decentralized control was posed as a CC selection problem.  Both theoretical and computational perspectives were provided, although the CC selection problem admits a degree of simplification in the discrete time setting (see Remark \ref{ShreyasWork} for details). The CC selection problem has been considered in \cite{Sezer} where a method for determining the minimum  number of essential inputs and outputs required for  decentralization was provided; however, the characterization does not cope with all cases, see, for instance \cite{Trave} (page 219).  Reference \cite{optimumfeedbackpattern} considers a CC selection problem with general heterogeneous communication costs subject to the constraint that the closed-loop system has no structurally fixed modes. The proposed solution is  suboptimal and, in particular, the framework  does not account for  actuator/sensor placement. In contrast, in this paper, we provide jointly optimal solutions for the I/O and CC selection problems assuming homogeneous actuator/sensor placement and communication costs. 

\vspace{-.4cm}

\subsection*{Main Contributions}

\vspace{-.1cm}

In this paper we propose  an unified framework to address  structural control system design, which  solves  the sparsest   I/O, as well as the jointly sparsest I/O and CC selection problems, by  exploiting the implications  of the former into the later. Moreover, the proposed solutions  are efficient because they can be implemented using polynomial algorithms in the number of the state variables.
The main results of this paper are outlined as follows: first, we provide a solution to (3a) and (3b) (by duality), which leads to the solution of $\mathcal P_1^d$. Next, inspired by the solution of $\mathcal P_1^d$ we provide a new characterization of structural controllability and structural observability, using the bipartite graph representation of the original system digraph and its directed acyclic graph representation. Then, by considering this new representation we compute and characterize $\bar B$ that solves (2a), and by duality obtain $\bar C$ that is a solution to (2b), and hence  the solution to $\mathcal P_1$. Next, we show that the solution to $\mathcal P_2$ can be obtained by using particular solutions to $\mathcal P_1$. Furthermore, in each of the above cases, we describe all possible solutions to the respective design problems. Finally, we provide polynomial complexity algorithmic procedures to compute solutions for the I/O and CC selection problems. Preliminary results concerning the solutions to $\mathcal{P}_{1}^{d}$, the dedicated I/O design problem, were presented in \cite{Pequito}. Here, in addition to addressing the design problems $\mathcal{P}_{1}$ and $\mathcal{P}_{2}$, we develop another equivalent characterization of solutions to $\mathcal{P}_{1}^{d}$ which yields simpler and with lower computational complexity algorithmic procedures to construct solutions to $\mathcal{P}_{1}^{d}$. Additionally, it acts as a bridge to more general constructs used in $\mathcal{P}_{1}$ and $\mathcal{P}_{2}$.

The rest of the  paper is organized as follows: Section \ref{NotationTerm} reviews some concepts and introduces fundamental  results in structural systems theory and establish their relations to graph-theoretic constructs. Subsequently, in Section \ref{secIII} we present a new necessary and sufficient condition for a system to be structurally controllable, and consequently structurally observable (by duality). These are used in Section  \ref{probp1} to describe  the solutions to problem $\mathcal{P}_1$ followed by Section \ref{probp2} where we provide  solutions to problem $\mathcal{P}_2$. In Section \ref{complexAnalysis} we present algorithmic procedures with polynomial complexity  to compute solutions to $\mathcal{P}_{1}^{d}$, $\mathcal{P}_{1}$ and $\mathcal{P}_{2}$.  Next, we present a detailed illustrative example in Section \ref{illustrativeexample} that explores the solutions to the different problems addressed in the paper. Finally,
Section \ref{secConclusion}  concludes the paper and discusses avenues for further research.
The proofs are relegated to the Appendix.

}

\vspace{-0.34cm}

\section{Preliminaries and Terminology}\label{NotationTerm}

\vspace{-0.2cm}

In this section we start to recall some classical concepts in structural systems \cite{dionSurvey} and some graph theoretic notions~\cite{Cormen:2001:IA:580470}. Consider  a linear time invariant (LTI) system
\vspace{-0.45cm}
\begin{equation}
\dot x= A x + Bu,\qquad\qquad y=Cx,\qquad\qquad   x(0)=x_0\in \mathbb{R}^{n}
\label{stateSpaceEq}
\vspace{-0.4cm}
\end{equation}
\noindent with $u\in\mathbb{R}^p$, $y\in\mathbb{R}^m$, and appropriate dimensions of the matrices $A,B,C$. In order to perform structural analysis efficiently, it is customary to associate \eqref{stateSpaceEq} with a directed graph (digraph) $\mathcal{D}=(\mathcal V,\mathcal E)$,  in which $\mathcal V$ denotes the set of \textit{vertices} and $\mathcal E$ represents the set of \textit{edges}, such that, an edge $(v_j,v_i)$ is directed from vertex $v_j$ to vertex $v_i$. To this end, let $\bar A \in \{0,1\}^{n\times n}$, $\bar B\in\{0,1\}^{n\times p}$ and $\bar C \in  \{0,1\}^{m\times n}$ be binary matrices that represent the structural
patterns of $A$, $B$ and $C$ respectively. Denote by $\mathcal{X}=\{x_1,\cdots,x_n\}$, $\mathcal{U}=\{u_1,\cdots,u_p\}$ and $\mathcal{Y}=\{y_1,\cdots,y_m\}$ the sets of state, input and output vertices, respectively. Denote by $ \mathcal{E}_{\mathcal{X},\mathcal{X}}=\{(x_i,x_j):\ \bar A_{ji}\neq 0\}$, $\mathcal{E}_{\mathcal{U},\mathcal{X}}=\{(u_j,x_i):\ \bar B_{ij}\neq 0\}$, $ \mathcal{E}_{\mathcal{X},\mathcal{Y}}=\{(x_i,y_j):\ \bar C_{ji}\neq 0\}$ and $ \mathcal{E}_{\mathcal{Y},\mathcal{U}}=\{(y_j,u_i):\ \bar K_{ij}\neq 0\}$ for a given \emph{information pattern} $\bar K\in \{0,1\}^{p\times m}$. We may then introduce the digraphs $\mathcal{D}(\bar A)=(\mathcal{X},\mathcal{E}_{\mathcal{X},\mathcal{X}})$, $\mathcal{D}(\bar A,\bar B)=(\mathcal{X}\cup \mathcal{U},\mathcal{E}_{\mathcal{X},\mathcal{X}}\cup \mathcal{E}_{\mathcal{U},\mathcal{X}} )$, $\mathcal{D}(\bar A,\bar C)=(\mathcal{X}\cup \mathcal{Y},\mathcal{E}_{\mathcal{X},\mathcal{X}}\cup \mathcal{E}_{\mathcal{X},\mathcal{Y}} )$, $\mathcal{D}(\bar A,\bar B,\bar C)=(\mathcal{X}\cup \mathcal{U}\cup \mathcal{Y},\mathcal{E}_{\mathcal{X},\mathcal{X}}\cup \mathcal{E}_{\mathcal{U},\mathcal{X}}\cup \mathcal{E}_{\mathcal{X},\mathcal{Y}} )$ and $\mathcal{D}(\bar A,\bar B,\bar K,\bar C)=(\mathcal{X}\cup \mathcal{U}\cup \mathcal{Y},\mathcal{E}_{\mathcal{X},\mathcal{X}}\cup \mathcal{E}_{\mathcal{U},\mathcal{X}}\cup \mathcal{E}_{\mathcal{X},\mathcal{Y}} \cup \mathcal{E}_{\mathcal{Y},\mathcal{U}})$. Note that in the digraph $\mathcal{D}(\bar{A},\bar{B})$, the input vertices representing the zero-columns of $\bar{B}$ correspond to isolated vertices. As such, the number of \emph{effective} inputs, i.e., the inputs which actually exert control, is equal to the number of non-zero columns of $\bar{B}$, or, in the digraph representation, the number of input vertices that are connected to at least one state vertex through an edge in $\mathcal{E}_{\mathcal{U},\mathcal{X}}$. A similar interpretation of \emph{effective} outputs holds. This has important implications while interpreting the solutions to problems $\mathcal{P}_{1}$ and $\mathcal{P}_{2}$: for instance, in $\mathcal{P}_{1}$, although the design matrices $\bar{B}$ and $\bar{C}$ are specified to be of size $n\times n$ for notational and technical convenience, in most practical cases the optimal design matrices (as characterized later) will turn out to be sparse with several zero-columns making the effective number of inputs and outputs much smaller than $n$. A digraph $\mathcal{D}_s=(\mathcal V_s,\mathcal E_s)$ with $\mathcal V_s\subset \mathcal V$ and $\mathcal E_s\subset \mathcal E$ is called a \textit{subgraph} of $\mathcal{D}$, and it is a \emph{strict subgraph} if $\mathcal V_s\subsetneq \mathcal V$ and/or $\mathcal E_s\subsetneq \mathcal E$. If $\mathcal V_s=\mathcal V$, $\mathcal{D}_s$ is said to \textit{span} $\mathcal{D}$. Finally, a subgraph with some property $P$ is \emph{maximal} if there is no other subgraph $\mathcal{D}_{s'}=(\mathcal V_{s'},\mathcal E_{s'})$ of $\mathcal{D}$, such that $\mathcal{D}_s$ is a strict subgraph of $\mathcal{D}_{s'}$ and $\mathcal{D}_{s^{\prime}}$ satisfies property $P$. A vertex in a graph with no incoming and outgoing edges is called an \emph{isolated vertex}. A sequence of directed edges $\{(v_1,v_2),(v_2,v_3),\cdots,(v_{k-1},v_k)\}$, in which all the vertices are distinct, is called \textit{an  elementary path} from $v_1$ to $v_k$. A vertex with an edge to itself (i.e., a \emph{self-loop}), or an elementary path from $v_1$ to $v_k$ comprising an additional edge $(v_k,v_1)$, is called a \emph{cycle}. A vertex $v$ is said to be \emph{reachable} from another vertex $w$ if there exists an elementary path from $w$ to $v$; in this case we say that $w$ \emph{reaches} $v$.

In addition, we will require the following graph theoretical notions~\cite{Cormen:2001:IA:580470}: A digraph $\mathcal{D}$ is said to be strongly connected if there exists an elementary path between any  pair of vertices. A \emph{strongly connected component} (SCC) is a maximal subgraph $\mathcal{D}_S=(\mathcal V_S,\mathcal E_S)$ of $\mathcal{D}$ such that for every $v,w \in \mathcal V_S$ there exists a path from $v$ to $w$ and from $w$ to $v$.  Visualizing each SCC as a virtual node (or supernode), one may generate a \textit{directed acyclic graph} (DAG), in which each super node corresponds to a single SCC and a directed edge exists between two SCCs \emph{if and only if} there exists a directed edge connecting vertices in the SCCs in the original digraph. It may be readily seen that the resulting DAG is acyclic, i.e., devoid of cycles. The DAG associated with a given $\mathcal{D}=(\mathcal V,\mathcal E)$ may be efficiently generated with complexity $\mathcal{O}(|\mathcal V|+|\mathcal E|)$~\cite{Cormen:2001:IA:580470}, where  $|\mathcal V|$ and $|\mathcal E|$  denote the number of vertices in $\mathcal V$ and the number of edges in $\mathcal E$ respectively.  The SCCs in a DAG may be classified as follows:

\vspace{-.2cm}

\begin{definition}[Non-top linked SCC/Non-bottom linked SCC]\label{linkedSCC}
An SCC is said to be linked  if it has at least one incoming/outgoing edge into its vertices to/from the vertices in another SCC. In particular,  an SCC is  \emph{non-top linked} if it has no incoming edge to its vertices from  another SCC and \emph{non-bottom linked} if it has no outgoing edge to another SCC.
\hfill $\diamond$
\end{definition}
We say that an SCC $\mathcal N=(\mathcal V,\mathcal E)$ is \emph{reachable} from another SCC $\mathcal N'=(\mathcal V',\mathcal E')$  if there exists a elementary path starting in $w\in \mathcal V'$ and ending in $v\in \mathcal V$; or,  we  say that the SCC $\mathcal N'$  \emph{reaches}  the SCC $\mathcal N$.

For any two vertex sets $\mathcal S_{1}, \mathcal S_{2}\subset \mathcal V$ we define the   \textit{bipartite graph} $\mathcal{B}(\mathcal S_1,\mathcal S_2,\mathcal E_{\mathcal S_1,\mathcal S_2})$ whose vertex set is given by $\mathcal S_{1}\cup \mathcal S_{2}$ and the edge set $\mathcal E_{\mathcal S_1,\mathcal S_2}\subset\{(s_1,s_2) \ :\ s_1 \in \mathcal S_1, s_2 \in \mathcal S_2  \ \}$. The bipartite graph $\mathcal B(\mathcal V,\mathcal V,\mathcal E)$ is said to be the bipartite graph associated with $\mathcal D(\mathcal V,\mathcal E)$.

Given $\mathcal{B}(\mathcal S_1,\mathcal S_2,\mathcal E_{\mathcal S_1,\mathcal S_2})$, a matching $M$ corresponds to a subset of edges in $\mathcal E_{\mathcal S_1,\mathcal S_2}$ that do not share vertices, i.e., given edges  $e=(s_1,s_2)$ and $e'=(s_1',s_2')$ with $s_1,s_1' \in \mathcal S_1$ and $s_2,s_2'\in \mathcal S_2$, $e, e' \in M$ only if $s_1\neq s_1'$ and $s_2\neq s_2'$. A maximum matching $M^{\ast}$ is  a matching $M$ that has the largest number of edges among all possible matchings. Given a bipartite graph $\mathcal{B}(\mathcal S_1,\mathcal S_2,\mathcal E_{\mathcal S_1,\mathcal S_2})$, the maximum matching problem may be solved  efficiently in $\mathcal{O}(\sqrt{|\mathcal S_1\cup \mathcal S_2|}|\mathcal E_{\mathcal S_1,\mathcal S_2}|)$ \cite{Cormen:2001:IA:580470}. Vertices in $\mathcal S_1$ and $\mathcal S_2$ are \textit{matched vertices} if they belong to an edge in the  matching $M$, otherwise, we designate them as \textit{unmatched vertices}. If there are no unmatched vertices, we say that we have a \textit{perfect match}. Notice  that a  matching $M$ (in particular a maximum matching) may not be unique and the above notions of matched or unmatched vertices are specific to a given matching.

For ease of referencing, in the sequel, the term \emph{right-unmatched vertices} (or \emph{left-unmatched vertices}) w.r.t. a  matching $M$ associated with  $\mathcal{B}(\mathcal S_1,\mathcal S_2,\mathcal E_{\mathcal S_1,\mathcal S_2})$ will refer to only those vertices in $\mathcal S_{2}$ (or $\mathcal S_1$) that do not belong to a matched edge in $M$.  Additionally, if we associate a \emph{weight} to each edge in the bipartite graph we may be interested in determine the maximum matching which the sum of the weights incurs in the minimum/maximum cost. These problems are known as the minimum/maximum weight maximum matching and are solved using, for instance, the Hungarian algorithm with computational complexity  $\mathcal O(\max \{\mathcal S_1,\mathcal S_2\}^3)$, see \cite{Cormen:2001:IA:580470} for details.

%
%
%
%

\vspace{-0.6cm}

\subsection{Structural Controllability and Observability }
\vspace{-0.2cm}
{
Given digraphs $\mathcal{D}(\bar A)$, $\mathcal{D}(\bar A,\bar B)$, $\mathcal{D}(\bar A,\bar C)$ or  $\mathcal{D}(\bar A,\bar B,\bar C)$ (when appropriate), we further define the following special subgraphs  \cite{Lin_1974}:

\noindent $\bullet$ \emph{State Stem} -  An isolated vertex or an elementary path, composed exclusively of state vertices.

\noindent $\bullet$ \emph{Input Stem} -  An input vertex linked to the root of a state stem.

\noindent $\bullet$ \emph{Output Stem} - A state stem linked from the tip to an output vertex.

\noindent $\bullet$ \emph{Input-Output Stem} -  An input vertex linked to the root of a state stem and linked from the state stem tip to an output vertex.

\noindent $\bullet$ \emph{State Cactus} -  Defined recursively  as follows: A state stem is a state cactus.  A state cactus connected  to a cycle from any state vertex  is also a state cactus.

\noindent $\bullet$ \emph{Input Cactus} - Defined recursively  as follows: An input stem with at least one state vertex is an input cactus. An input cactus connected  to a cycle (comprised of state vertices only) from any vertex (either state or input vertex) is also an input cactus.

\noindent $\bullet$ \emph{Output Cactus} - Defined recursively  as follows: An output stem with at least one state vertex is an output cactus. A cycle connected to an output cactus at any vertex (either state or output vertex) is also an output cactus.

\noindent $\bullet$ \emph{Input-Output Cactus} - Defined recursively  as follows: An input-output stem with at least one state vertex is an input-output cactus. An input-output cactus connected  to/from a cycle from/to any of its vertices is also an input-output cactus.

\noindent $\bullet$ \emph{Chain} - A single cycle or a  group of disjoint cycles (composed of state vertices) connected to each other in a sequence. In other words, a DAG where each supernode is a cycle.

The root and the tip of a stem are also the root and the tip of the associated cactus. Note that, by definition, an input cactus may have an input vertex linked to several state vertices, which means,  for example, that the input vertex may connect to the root of a state stem and also be linked to one or more states in a chain.

A system $(\bar A,\bar B)$ is said to be structurally controllable if there exists a pair   $(A_0,B_0)$ of real matrices   with the same structure, i.e., location of zero and non-zero entries, as $(\bar A,\bar B)$ such that $(A_0,B_0)$ is controllable \cite{dionSurvey}. Structural controllability may be characterized as follows:

\vspace{-0.2cm}

\begin{theorem}[\hspace{-0.02cm}\cite{dionSurvey}]\label{Theorem1}
For LTI systems described by \eqref{stateSpaceEq}, the following statements are equivalent:
\begin{itemize}
\item[i)] The corresponding structured linear system $(\bar A,\bar B)$ is structurally controllable.
\item[ii)] The digraph $\mathcal{D}(\bar A,\bar B)$ is spanned by a disjoint union of input cacti.
\hfill$\diamond$
\end{itemize}
\end{theorem}

}
{\setstretch{1.47}
Similarly, $ (\bar A, \bar C)$ is structurally observable if and only if $(\bar A^T,\bar C^T)$ is structurally controllable, which is equivalent to $\mathcal{D}(\bar{A},\bar{C})$ being spanned by a disjoint union of output cacti.

%

\vspace{-0.4cm}

\subsection{ Structural Fixed Modes}\label{decControlWang}

Let \hspace{0.1cm} $
[\bar M]_{n_1,n_2}=\left\{M\in \mathbb{R}^{n_1\times n_2} : \ M_{ij}= 0\  \text{if} \  \bar M_{ij}=0, \ 0 \le i \le n_1, 0\le j \le n_2 \right\}$ \hspace{0.1cm} 
denote\hspace{0.1cm}  the \\ equivalence class of matrices with a given structure $\bar{M}$ of dimensions $n_1\times n_2$, and $\sigma(X)$ denote the set of eigenvalues of a square matrix $X$.
\vspace{-0.3cm}

\begin{definition}[\hspace{-0.1cm} \cite{Davison}] Consider system  \eqref{stateSpaceEq} with p inputs and m outputs.
The set
$
\sigma_{\bar K}=\bigcap_{K \in [ \bar K]_{p,m}}  \sigma(A+BKC)
$
is defined to be the set of fixed modes of  the closed-loop of \eqref{stateSpaceEq} w.r.t. the information pattern $\bar K$, where $[\bar K]_{p,m}$ is the set of all possible $p \times m$ constant output feedback matrices.
\hfill $\diamond$
\end{definition}
\vspace{-0.3cm}

The stabilizability of a system  under an information pattern $\bar K$ is related to the fact that there are no fixed modes in $\sigma_{\bar K}$ with nonnegative real part.  In particular, in \cite{Davison} it is shown that the condition $\sigma_{\bar K} =\emptyset$  is both necessary and sufficient for \emph{almost arbitrary} pole placement with output feedback.

Fixed modes also have a structural counterpart, the \textit{structurally fixed modes} (SFM), given next.
\vspace{-0.3cm}
\begin{definition}[\hspace{-0.02cm}\cite{Papadimitriou84asimple}]
System \eqref{stateSpaceEq} in closed-loop, denoted by $(\bar{A},\bar{B},\bar{K},\bar{C})$, is said to have \emph{structurally fixed modes} w.r.t. an information pattern $\bar K$ if for all $A \in [\bar A]_{n\times n}$, $B \in [\bar B]_{n\times p}$, $C \in [\bar C]_{m\times n}$, we have
$
\bar\sigma_{\bar K}=\bigcap_{K \in [\bar K]_{p\times m}}  \sigma(A+BKC) \neq \emptyset .
$\hfill $\diamond$
\end{definition}

\vspace{-0.3cm}

Conversely, a structural system $(\bar{A},\bar{B},\bar K, \bar{C})$ has no structurally fixed modes, if there exists at least one instantiation $A\in [\bar{A}]$, $B\in [\bar{B}]$, $C\in [\bar{C}]$ which has no fixed modes, i.e., $\cap_{K\in [\bar{K}]}\sigma(A+BKC)=\emptyset$. In this latter case, it may be shown (see \cite{sezerFixedModes}) that almost all systems in the sparsity class $(\bar{A},\bar{B},\bar{C})$ have no fixed modes with respect to $\bar{K}$, and, hence, allow pole-placement arbitrarily close to any pre-specified set of eigenvalues. This is the key motivation behind our constraint of designing systems with no SFMs in the CC selection problem in $\mathcal P_2$.

Now, consider the following graph-theoretic conditions that ensure the absence of structrurally fixed modes.

\vspace{-0.3cm}

\begin{theorem}[\hspace{-0.02cm}\cite{Pichai3}]
The structural system $(\bar A,\bar B,\bar C)$ associated with \eqref{stateSpaceEq} has no structurally fixed modes w.r.t. an information pattern $\bar K$, \underline{if and only if} both of the following conditions hold:
\begin{enumerate}
\item[a)] in $\mathcal D(\bar A,\bar B,\bar K,\bar C)=(\mathcal X\cup\mathcal U\cup \mathcal Y,\mathcal E_{\mathcal X,\mathcal X}\cup\mathcal E_{\mathcal X,\mathcal Y}\cup\mathcal E_{\mathcal U,\mathcal X}\cup\mathcal E_{\mathcal Y,\mathcal U})$, each state vertex $x\in \mathcal{X}$ is contained in an SCC which includes an edge of $\mathcal{E}_{\mathcal{Y},\mathcal{U}}$;
\item[b)] there exists a finite disjoint union of cycles $\mathcal{C}_k=(\mathcal{V}_k,\mathcal{E}_k)$ in $\mathcal{D}(\bar A,\bar B,\bar K,\bar C)$ with $k\in\mathbb{N}$ such that
$
\mathcal{X}\subset \bigcup_{k} \ \mathcal V_k\ .
$\hfill $\diamond$
\end{enumerate}
\label{Theorem2}
\end{theorem}
}

\vspace{-1cm}

\section{Solution to Problem $\mathcal{P}_1^d$}\label{secIII}
\vspace{-0.1cm}

In this section we introduce a new characterization of structural controllability/observability by considering the bipartite representation of the system digraph $\mathcal D(\bar A)$ and its directed acyclic graph (DAG) representation. We start by describing a result that provides a bridge between structural systems concepts and graph theoretic constructs, such as the maximum matching problem. This is fundamental to provide and explicitly  characterize  solutions to $\mathcal P_1^d$, and consequently  solutions to both $\mathcal P_1$ and $\mathcal P_2$. The proofs are based on standard graph theoretic properties and  relegated to the Appendix.


We start by a couple of useful results about general digraphs properties.
\vspace{-0.3cm}

\begin{lemma}\label{lemmaA}
Let $\mathcal D=(\mathcal X,\mathcal E_{\mathcal X,\mathcal X})$ be an SCC which  is spanned by a disjoint union of cycles $\mathcal C^i$, $i=1,\ldots,c_{\mathcal D}$. Then, for each such cycle $\mathcal{C}^{i}$ there exists a chain that spans $\mathcal{D}$ and whose first element is $\mathcal{C}^{i}$. \hfill $\diamond$
\end{lemma}

\vspace{-0.4cm}

The following extension of Lemma 1 to arbitrary digraphs holds.

\vspace{-0.3cm}

\begin{lemma}\label{lemmaB}
Let $\mathcal D=(\mathcal X,\mathcal E_{\mathcal X,\mathcal X})$ be a digraph  comprising $\gamma$ non-top linked SCCs $\mathcal N^T_i$ with $i=1,\ldots,\gamma$. Further, let $\mathcal{D}$ be spanned by a disjoint collection of cycles, and, $\{ \mathcal C_i'\}_{i=1,\ldots,\gamma} $ be any subcollection of such cycles such that each $\mathcal{C}^{\prime}_{i}\in\mathcal{N}_{i}^{T}$, $i=1,\cdots,\gamma$. Then, there exists a collection of disjoint $\gamma$ chains, such that the $i$-th chain has its first element as the cycle $\mathcal{C}^{\prime}_{i}$ for $i=1,\cdots,\gamma$. \hfill $\diamond$
\end{lemma} 

\vspace{-0.3cm}

In addition, consider the following result that relates a maximum matching of the state biparte representation with the state digraph.

\vspace{-0.3cm}

\begin{lemma}[Maximum Matching Decomposition]\label{Lemma2}
Consider the digraph $\mathcal{D}(\bar A)=(\mathcal{X},\mathcal{E}_{\mathcal{X},\mathcal{X}})$ and let $M^*$ be a maximum matching associated with the bipartite graph $\mathcal{B}(\mathcal{X},\mathcal{X},\mathcal{E}_{\mathcal{X},\mathcal{X}})$. Then, the digraph $\mathcal D^*=(\mathcal X, M^*)$  constitutes a disjoint union of cycles and  state stems (by convenction an isolated state vertex is a state stem), with roots in the right-unmatched vertices and tips in the left-unmatched vertices of $M^{\ast}$, that  span $\mathcal{D}(\bar A)$. Moreover, such a decomposition is \emph{minimal}, in the sense that, no other spanning subgraph decomposition of $\mathcal{D}(\bar{A})$ into state stems and cycles contains strictly fewer number of state stems.\hfill $\diamond$
\end{lemma}
\vspace{-0.3cm}

{\setstretch{1.44}
Remark that as a consequence of Lemma \ref{Lemma2}, if the set of unmatched vertices is empty, then the original graph is spanned by a disjoint union of cycles comprising the edges in a maximum matching. This result  coincides with a result previously established by K\"{o}nig  (see the appendix in \cite{Reinschke:1988}). In addition, note that Lemma \ref{Lemma2} states that the maximum matching problem leads to two different kinds of matched edge sequences in $M^*$; sequences of edges in $M^*$ starting in right-unmatched state vertices (i.e., state stems with more than one state vertex),  and the remaining  sequences of edges that start and end in matched vertices  (i.e., cycles).  Recall that the maximum matching is not unique, which implies that the above decomposition into state stems and cycles is not unique as well. Note that, Lemma \ref{Lemma2} provides no information on the length of each state stem or the number of cycles that constitute the maximum matching decomposition.

 A \emph{feasible dedicated input configuration} is a subset $\mathcal S_u\subset \mathcal X$ of state variables to which by assigning  dedicated inputs we obtain a structurally controllable system. Note that each  feasible dedicated input configuration corresponds to a unique canonical $\bar{B}$ up to column permutation.  The next result provides a characterization of feasible dedicated input configurations.

\vspace{-0.3cm}

\begin{theorem}[Feasible dedicated input configuration]\label{FDIC}
Let $\mathcal D(\bar A)=(\mathcal X,\mathcal E_{\mathcal X,\mathcal X})$ denote the system digraph and denote by $\mathcal B \equiv \mathcal B(\mathcal X,\mathcal X,\mathcal E_{\mathcal X,\mathcal X})$ the  associated bipartite representation. A set $\mathcal S_u \subset \mathcal X$ is a feasible dedicated input configuration \underline{if and only if}   there exist subsets $\mathcal U_R$ and $\mathcal{A}_{u}$ of $\mathcal{S}_{u}$, such that $\mathcal{U}_{R}$ corresponds to the set of right-unmatched vertices of some maximum matching of $\mathcal B$, and  $\mathcal A_u$ contains one state variable from each non-top linked SCC of $\mathcal D(\bar A)$.\hfill $\diamond$
\end{theorem}

\vspace{-0.3cm}

By duality we obtain a similar characterization of a \emph{feasible dedicated output configuration},  which is defined to be a subset of state variables, to which  assigning  dedicated outputs leads to a structurally observable system.


From the proof of Theorem \ref{FDIC} (respectively the dual of Theorem \ref{FDIC}), it is easy to realize that the SCCs that contain a right-unmatched (respectively left-unmatched) vertex of some maximum matching of $\mathcal{B}$ play an important role in the selection of feasible dedicated input (respectively output) configurations. This motivates  the following notions.

\vspace{-0.3cm}

\begin{definition}\label{topassignedSCC}
Let $\mathcal{D}(\bar A)=(\mathcal{X},\mathcal{E}_{\mathcal{X},\mathcal{X}})$ and $M^*$ be a maximum matching associated with its bipartite representation.
A non-top linked SCC is said to be a \emph{top assignable SCC} with respect to $M^{\ast}$ if it contains at least one right-unmatched vertex in $M^{\ast}$. Similarly, a non-bottom linked SCC is said to be a \emph{bottom assignable} SCC  with respect to $M^{\ast}$ if it contains at least one left-unmatched vertex.
\hfill $\diamond$
\end{definition}

\vspace{-0.3cm}

Note that the total number of top and bottom assignable SCCs  may depend on the particular maximum matching $M^{\ast}$ (not unique in general) under consideration; as such we introduce the following definition:

\begin{definition}\label{maxtopassign}
Consider the digraph $\mathcal{D}(\bar A)=(\mathcal{X},\mathcal{E}_{\mathcal{X},\mathcal{X}})$.
The \emph{maximum top assignability  index} of $\mathcal{D}(\bar A)$ is the maximum number of top assignable SCCs among the  maximum matchings $M^*$ associated with $\mathcal{B}(\mathcal{X},\mathcal{X},\mathcal{E}_{\mathcal{X},\mathcal{X}})$.  Similarly, the \emph{maximum bottom assignability  index} of $\mathcal{D}(\bar A)$ is the maximum number of bottom assignable SCCs among the  maximum matchings $M^*$ associated with $\mathcal{B}(\mathcal{X},\mathcal{X},\mathcal{E}_{\mathcal{X},\mathcal{X}})$.
\hfill $\diamond$
\end{definition}
}
\vspace{-0.3cm}



A maximum matching that attains the top assignability (respectively bottom assignability) index is said to be top (respectively bottom) assignable. The following result  addresses (3a)  by providing the minimum number  of required  dedicated inputs, hence the minimum number of columns in $\bar B$ (each with only one non-zero entry) required to ensure the structural controllability of the pair $(\bar{A}, \bar{B})$.

\vspace{-0.4cm}

\begin{theorem}[Minimum number of dedicated inputs]
Let $\mathcal{D}(\bar A)=(\mathcal{X},\mathcal{E}_{\mathcal{X},\mathcal{X}})$ be the system digraph with $\beta$ denoting the number of non-top linked SCCs in its DAG representation.
Then, the minimum number of dedicated inputs $p$ required to achieve structural controllability is given by
\vspace{-0.6cm}
\begin{equation}
p=m+\beta-\alpha,
\vspace{-0.4cm}
\label{eq2}
\end{equation}
where $m$ denotes the number of right-unmatched vertices in any maximum matching of the bipartite graph $\mathcal B(\mathcal X,\mathcal X,\mathcal E_{\mathcal X,\mathcal X})$ and  $\alpha$ denotes the maximum top assignability  index of $\mathcal{D}(\bar A)$.
\hfill $\diamond$
\label{Theorem4}
\end{theorem}
\vspace{-0.2cm}

It may be readily verified from the definitions, that if $\mathcal{D}(\bar A)$ is strongly connected,  we have $\beta=1$, and $\alpha$ may only assume two values, $0$ or $1$, depending on whether $m=0$ or $m\neq 0$, respectively. As such,  Theorem \ref{Theorem4} may be simplified significantly if $\mathcal{D}(\bar A)$ is known to be strongly connected, in which case $p=\max(m,1)$ and  coincides with the result presented in \cite{liu11}.

In addition, Theorem \ref{FDIC} and Theorem \ref{Theorem4} provide us with the following characterization of a \emph{minimal feasible dedicated input configuration}, i.e., the subset with the minimum number of state variables to which we need to assign dedicated inputs to ensure structural controllability.

\vspace{-0.2cm}

\begin{theorem}[Minimal feasible dedicated input configuration]\label{MFDIC}
Let $\mathcal D(\bar A)=(\mathcal X,\mathcal E_{\mathcal X,\mathcal X})$ denote the system digraph and $\mathcal B \equiv \mathcal B(\mathcal X,\mathcal X,\mathcal E_{\mathcal X,\mathcal X})$ its bipartite representation. A set $\mathcal S_u\subset \mathcal X$ is a minimal feasible dedicated input configuration \underline{if and only if} there exist two disjoint  subsets $\mathcal U_R$ and $\mathcal{A}_{u}^c$ such that $\mathcal{S}_{u}=\mathcal{U}_{R}\cup\mathcal{A}_{u}^c$, $\mathcal{U}_{R}$ corresponds to the set of right-unmatched vertices of some maximum matching of $\mathcal B$ with maximum assignability, and  $\mathcal A_u^c$ comprises of only one state variable from each non-top linked SCC of $\mathcal D(\bar A)$ not assigned by the right-unmatched vertices in $\mathcal U_R$.
\hfill $\diamond$
\end{theorem}
\vspace{-0.2cm}

Theorem \ref{MFDIC}  will be used  to describe all possible solutions of (3a), and by duality of linear systems,  it can  also be used to describe the solutions of (3b). Therefore,  the following result  holds.
\vspace{-0.4cm}

\begin{theorem}\label{solP1d}
A structural pair $(\bar{B},\bar{C})$ is a solution to $\mathcal{P}_{1}^{d}$ \underline{if and only if} $\bar{B}$ and $\bar{C}$ correspond to the dedicated assignment in Theorem~\ref{MFDIC} and its dual result, respectively. \hfill $\diamond$
\end{theorem}
\vspace{-0.3cm}

A polynomial complexity algorithm for explicitly constructing minimal feasible dedicated configurations (and hence structural pairs $(\bar{B},\bar{C})$ solving $\mathcal{P}_{1}^{d}$) is provided in Section VI.  In the next section we show how  minimal feasible dedicated input  (and output) configuration can be used to obtain solutions $(\bar B,\bar C)$  to $\mathcal P_1$.

\vspace{-0.4cm}

\section{Solution to Problem $\mathcal{P}_1$}\label{probp1}

Intuitively, the inputs in a minimum feasible dedicated input configuration may be classified into the following two types: those dedicated inputs that may be merged into a common input such that structural controllability is retained, and those that cannot. Based on the above classification, the solution to problem (2a) is achieved in Theorem~\ref{sol2a}. Subsequently, as a corollary to Theorem~\ref{sol2a} (Corollary \ref{minsol2a}) we further obtain the \emph{minimal} solution to (2a), i.e., the solution of (2a) with the minimum number of  non-zero columns of $\bar B$, in other words, the minimum number of inputs that actuate at least one state variable. By duality, the results obtained above for the input design problem (2a) and its variants are extended to the output design (2b) and its variants.
Finally, by using Theorem \ref{sol2a} and its dual we obtain the solution to $\mathcal{P}_{1}$ in Theorem \ref{solP1}.

We start this section by providing a new characterization of structural controllability in terms of the system bipartite representation and associated maximum matchings.

\vspace{-0.2cm}

\begin{theorem}\label{structuralControl}
Let $\mathcal D(\bar A,\bar B)=(\mathcal X\cup \mathcal U,\mathcal E_{\mathcal X,\mathcal X}\cup \mathcal E_{\mathcal U,\mathcal X})$ denote the state-input digraph and $\mathcal B \equiv \mathcal B(\mathcal X,\mathcal X,\mathcal E_{\mathcal X,\mathcal X})$ the state bipartite representation. The pair $(\bar A,\bar B)$ is structurally controllable \underline{if and only if} there exist the following three  subsets: $\mathcal U_R\subset \mathcal X$ corresponding to the set of right-unmatched vertices of some maximum matching of $\mathcal B$; $\mathcal A_u\subset \mathcal  X$ comprising one state variable from each non-top linked SCC of $\mathcal D(\bar A)$; and $\mathcal U_S\subset \mathcal U$ such that:
\begin{itemize}
\item[(i)] to each $x\in\mathcal{U}_{R}$ there exists a distinct $u\in\mathcal{U}_{S}$ assigned to $x$, i.e., with $(u,x)\in\mathcal{E}_{\mathcal{U},\mathcal{X}}$; and
\item[(ii)] to each $x\in\mathcal{A}_{u}$ there exists a $u\in\mathcal{U}_{S}$ with  $(u,x)\in\mathcal{E}_{\mathcal{U},\mathcal{X}}$.\hfill $\diamond$
\end{itemize}
\end{theorem}

\vspace{-0.5cm}
{

\begin{remark}
\label{rem:struc}
Note that, by \emph{distinct} in condition (i) of Theorem~\ref{structuralControl}, we mean that if $x$ and $\acute{x}$ are two distinct states in $\mathcal{U}_{R}$, there exist $u$ and $\acute{u}$ in $\mathcal{U}_{S}$ with $u\neq\acute{u}$ such that the edges $(u,x)$ and $(\acute{u},\acute{x})$ are in $\mathcal{E}_{\mathcal{U},\mathcal{X}}$. However, no such distinction is required in the input assignment to states in $\mathcal{A}_{u}$; in particular, connecting the inputs assigned to states in $\mathcal{U}_{R}$ to states in $\mathcal{A}_{u}\setminus\mathcal{U}_{R}$ is allowable as far as condition (ii) of Theorem~\ref{structuralControl} is concerned.

It follows from the above discussion that for a structurally controllable pair $(\bar{A},\bar{B})$, the number of effective inputs is at least $m$ if $m$, the number of right unmatched vertices in a maximum matching of $\mathcal{B}$, is non-zero, or, it is at least one otherwise; in particular, we have the lower bound $\max(m,1)$ on the number of effective inputs if $(\bar{A},\bar{B})$ is structurally controllable. \hfill $\diamond$
\end{remark}

\vspace{-0.2cm}

%

We now introduce some additional notation. Recall that an input structural matrix $\bar{B}\in\{0,1\}^{n\times n}$ may be equivalently specified by the edge set $\mathcal{E}_{\mathcal{U},\mathcal{X}}$ in the digraph representation $\mathcal{D}(\bar{A},\bar{B})$. In the following we will use the notation $\mathcal{E}_{\mathcal{U},\mathcal{X}}(\bar{B})$ to make this connection explicit. For a maximum matching $M$ of the state bipartite representation $\mathcal{B}$, let $\mathcal{U}_{R}(M)=\{x_{i_{1}}(M),\cdots,x_{i_{m}}(M)\}$ be an enumeration of the state variables corresponding to its $m$ right-unmatched vertices and let $\{\bar{\mathcal{N}}^{T}_{i_{k}}(M)\}_{k=1}^{\gamma(M)}$ be the collection of non-assigned  non-top linked SCCs of $\mathcal{D}(\bar{A})$, i.e.,  such that $\bar{\mathcal{N}}^{T}_{i_{k}}(M)\cap\mathcal{U}_{R}(M)=\emptyset$ for all $k=1,\cdots,\gamma(M)$; clearly, the number of non-assigned non-top linked SCCs $\gamma(M)\leq\beta$, where $\beta$ denotes the total number of non-top linked SCCs of $\mathcal{D}(\bar{A})$. For each such maximum matching $M$, define the set $\mathcal{I}(M)$ of input structural matrices $\bar{B}\in\{0,1\}^{n\times n}$ representing the edges from different inputs to the right-unmatched vertices associated with $M$ and an edge from an input to a single state variable in each  non-assigned non-top linked SCC, as
\vspace{-0.3cm}
\begin{align}
\label{def:I_M} \mathcal{I}(M)=&\left\{\bar{B}~:~\mathcal{E}_{\mathcal{U},\mathcal{X}}(\bar{B})=\left\{\cup_{j=1}^{m}\left(u_{i_{j}},x_{i_{j}}(M)\right)\right\}\cup\left\{\cup_{k=1}^{\gamma(M)}\left(\acute{u}_{i_{k}},\acute{x}_{i_{k}}\right)\right\} \right. \\
& \quad\nonumber \left. \mbox{s.t. $u_{i_{j}},\acute{u}_{i_{k}}\in\mathcal{U}$ for all $j=1,\cdots,M$, \ ~$k=1,\cdots,\gamma(M)$ and $u_{i_{j}}\neq u_{i_{\acute{j}}}$ for $j\neq\acute{j}$,} \right. \\ 
&\quad\quad\qquad\qquad\qquad\qquad\qquad\qquad\qquad\qquad \nonumber \left. \mbox{$\acute{x}_{i_{k}}\in\bar{\mathcal{N}}^{T}_{i_{k}}(M)$ for all $k=1,\cdots,\gamma(M)$}\right\}.
\end{align}

\vspace{-0.6cm}

\begin{remark}\label{remarkStruct}
Using the above construction, Theorem~\ref{structuralControl} may be restated as follows: the pair $(\bar{A},\bar{B})$ is structurally controllable if and only if there exist a maximum matching $M$ of $\mathcal{B}$ and a structural input matrix $\acute{\bar{B}}\in\mathcal{I}(M)$ such that $\acute{\bar{B}}\leq\bar{B}$ (where the inequality is to be interpreted entry-wise). 
\vspace{-0.8cm}

\hfill $\diamond$
\end{remark}

\vspace{-0.7cm}

Now, note that, from~\eqref{def:I_M} we have for any $\acute{\bar{B}}\in\mathcal{I}(M)$,
\vspace{-0.4cm}
\begin{equation}
\label{def:I_M1}\|\acute{\bar{B}}\|_{0}=m+\gamma(M),
\vspace{-0.4cm}
\end{equation}
\noindent which, by the  characterization of structural controllability stated in Remark \ref{remarkStruct}, yields that

\vspace{-0.8cm}

\begin{equation}
\label{def:I_M2}\|\bar{B}\|_{0}\geq m+\min_{M}\gamma(M),
\vspace{-0.3cm}
\end{equation}
for any $\bar{B}$ such that $(\bar{A},\bar{B})$ is structurally controllable. Further, noting that  $\min_{M}\gamma(M)=\beta-\alpha$ (where $\alpha$ denotes the maximum top assignability index of $\mathcal{D}(\bar{A})$, see Definition~\ref{maxtopassign}) and the minimizers are those maximum matchings which have maximum top assignability, we immediately obtain the solution to (2a) as follows.
}

\vspace{-0.3cm}

\begin{theorem}[Solution to (2a)]\label{sol2a}
A structural matrix $\bar{B}$ is a solution to (2a) if and only if there exists a maximum matching $M$ of $\mathcal{B}$ with maximum top assignability such that $\bar{B}\in\mathcal{I}(M)$, where $\mathcal{I}(M)$ is defined in~\eqref{def:I_M}.
\hfill $\diamond$
\end{theorem}
\vspace{-0.4cm}


To obtain minimal solutions of (2a), i.e., the ones with the smallest number of effective inputs, let us define for each maximum matching $M$ of $\mathcal{B}$ the (non-empty) subset $\acute{\mathcal{I}}(M)$ of $\mathcal{I}(M)$ as follows:

If $m\neq 0$, let
\vspace{-.3cm}
\begin{align}
\label{def:I_Macute} \acute{\mathcal{I}}(M)=&\left\{\bar{B}~:~\mathcal{E}_{\mathcal{U},\mathcal{X}}(\bar{B})=\left\{\cup_{j=1}^{m}\left(u_{i_{j}},x_{i_{j}}(M)\right)\right\}\cup\left\{\cup_{k=1}^{\gamma(M)}\left(\acute{u}_{i_{k}},\acute{x}_{i_{k}}\right)\right\} \right. \\ \nonumber&\qquad\qquad \left. \mbox{s.t. $u_{i_{j}}\in\mathcal{U}$ for all $j=1,\cdots,M$ and $u_{i_{j}}\neq u_{i_{\acute{j}}}$ for $j\neq\acute{j}$,} \right. \\ \nonumber &\qquad\qquad\left. \mbox{$\acute{u}_{i_{k}}\in\{u_{i_{1}},u_{i_{2}},\cdots,u_{i_{M}}\}$ and $\acute{x}_{i_{k}}\in\bar{\mathcal{N}}^{T}_{i_{k}}(M)$  for all $k=1,\cdots,\gamma(M)$}\right\},
\end{align}
\vspace{-0.9cm}
or, if $m=0$, let
\vspace{-0.2cm}
\begin{align}
\label{def:I_Macute1}\qquad\qquad\qquad\qquad \acute{\mathcal{I}}(M)=&\left\{\bar{B}~:~\mathcal{E}_{\mathcal{U},\mathcal{X}}(\bar{B})=\left\{\cup_{k=1}^{\gamma(M)}\left(\acute{u},\acute{x}_{i_{k}}\right)\right\} \right. \\ \nonumber &\qquad\qquad \left. \mbox{s.t. $\acute{u}\in\mathcal{U}$ and $\acute{x}_{i_{k}}\in\bar{\mathcal{N}}^{T}_{i_{k}}(M)$ for all $k=1,\cdots,\gamma(M)$}\right\}.
\end{align}
\vspace{-1cm}

Clearly, for each maximum matching $M$, the set $\acute{\mathcal{I}}(M)$ defined in~\eqref{def:I_Macute}-\eqref{def:I_Macute1} (depending on whether $m\neq 0$ or not) is non-empty and coincides with the subset of those structural input matrices in $\mathcal{I}(M)$ with exactly $\max(m,1)$ number of effective inputs. Hence, by the lower bound in Remark~\ref{rem:struc} and the characterization in Theorem~\ref{sol2a}, we obtain the minimal solution to (2a) as follows:

\vspace{-0.3cm}

\begin{corollary}[Minimal solution to (2a)]\label{minsol2a}
A structural matrix $\bar{B}$ is a minimal solution to (2a) if and only if there exists a maximum matching $M$ of $\mathcal{B}$ with maximum top assignability such that $\bar{B}\in\acute{\mathcal{I}}(M)$, where $\acute{\mathcal{I}}(M)$ is defined in~\eqref{def:I_Macute}-\eqref{def:I_Macute1}.
\hfill $\diamond$
\end{corollary}
\vspace{-0.2cm}

In fact, in Corollary~\ref{minsol2a}, we provide a stronger result: we show that the set of structural input  matrices that solve $\mathcal P_1$ has a non-empty intersection with the set of structural input matrices  with the minimum number of non-zero columns (respectively rows) achieving structural controllability (respectively observability). Moreover, in the same corollary we explicitly characterize solutions of $\mathcal P_1$ that additionally possess the latter property, i.e., sparsest structural input matrices together with the minimum number of non-zero columns achieving structural controllability. By the above it readily follows, in particular, that the minimum number of effective inputs required to make a system structurally controllable is equal to $\max(m,1)$. This particular result was also obtained in \cite{dionSurveyKyb,Murota:2009:MMS:1822520}. However, in addition we obtain the sparsest design, i.e., we show that the minimum number of \emph{links} required between inputs and states to achieve structural controllability is $m+\beta-\alpha$ and explicitly characterize all such sparsest input configurations (see Theorem~\ref{sol2a}), which was not addressed in~\cite{dionSurveyKyb,Murota:2009:MMS:1822520}.

By duality in linear systems we can derive a solution to (2b), hence the solution to $\mathcal P_1$.
\vspace{-0.3cm}

\begin{theorem}[Solution to $\mathcal P_1$]\label{solP1}
The pair $(\bar B,\bar C)$ is a solution to $\mathcal P_1$ \underline{if and only if} $\bar B$ (resp. $\bar C$) is designed such that the characterization in Theorem \ref{sol2a} (resp. dual of Theorem \ref{sol2a}) holds. \hfill $\diamond$
\end{theorem}
\vspace{-0.3cm}

Using Corollary~\ref{minsol2a} we can derive the following result that plays a key role in obtaining the solution to $\mathcal P_2$.
\vspace{-0.3cm}

\begin{corollary}[Minimal solution to $\mathcal P_1$]\label{minsolP1}
The pair $(\bar B,\bar C)$ is a minimal solution to $\mathcal P_1$ \underline{if and only if} $\bar B$ (resp. $\bar C$) is designed such that the characterization in Corollary~\ref{minsol2a} (resp. dual of Corollary~\ref{minsol2a}) holds. \hfill $\diamond$
\end{corollary}
\vspace{-0.4cm}

In section \ref{complexAnalysis} we provide algorithmic procedures to efficiently (polynomial in the number of state variables) compute a minimal feasible dedicated input/output configuration, from which an efficient solution to problem $\mathcal{P}_1$ is obtained by using the characterization in Theorem \ref{solP1}.

\vspace{-0.4cm}

\section{Solution to Problem $\mathcal{P}_2$}\label{probp2}
{\setstretch{1.46}
Broadly, in this section, we will show that all possible solutions $(\bar{B},\bar{K},\bar{C})$ of $\mathcal{P}_{2}$ may be obtained (see Theorem \ref{Theorem8} for a formal statement) by considering minimal solutions $(\bar{B},\bar{C})$ of $\mathcal{P}_{1}$ and appropriately adding feedback edges between effective outputs in $\bar{C}$ and effective inputs in $\bar{B}$ through a procedure to be referred to as \emph{mix-pairing}. Given a solution to $\mathcal{P}_{1}$ such a mix-pairing characterizes the minimum number of feedback links that are required to ensure the requirements in Theorem \ref{Theorem2} which provides necessary and sufficient conditions for generic pole placement. To define such a mix-pairing procedure in full generality (so as to characterize all possible solutions of $\mathcal{P}_{2}$) and establish its minimality, we need several intermediate constructions and results detailed in the following.

As will be seen later, to achieve minimality in the feedback design procedure, the first step is to characterize a certain decomposition of the digraph $\mathcal{D}(\bar{A},\bar{B},\bar{C})$ into cycles and input-output stems, where $(\bar{B},\bar{C})$ is a solution to $\mathcal{P}_{1}
$. From a design point of view (and since we are interested in characterizing all solutions to $\mathcal{P}_{2}$), suppose $\bar{B}$ and $\bar{C}$ are obtained using two different maximum matchings of the state bipartite graph, say $M_{1}$ and $M_{2}$ respectively, as in Corollary \ref{minsolP1}. Clearly, $M_{1}$ (respectively $M_{2}$) provides a decomposition of the digraph $\mathcal{D}(\bar{A},\bar{B})$ (respectively $\mathcal{D}(\bar{A},\bar{C})$) into a disjoint union of cycles and input stems (respectively output stems); however, a decomposition of the joint digraph $\mathcal{D}(\bar{A},\bar{B},\bar{C})$ into a disjoint union of cycles and input-output stems may not be obvious (i.e., whether such a decomposition exists or its characterization) given the separate decompositions of ${D}(\bar{A},\bar{B})$ and $\mathcal{D}(\bar{A},\bar{C})$ (unless the maximum matchings $M_{1}$ and $M_{2}$ are equal). To this end, we provide a general graph-theoretic result which, given $M_{1}$ and $M_{2}$ (associated with $\bar{B}$ and $\bar
{C}$ respectively), characterizes a \emph{common} maximum matching $M^{\ast}$ of the state bipartite graph which explicitly provides a decomposition of the joint digraph $\mathcal{D}(\bar{A},\bar{B},\bar{C})$  into cycles and input-output stems. The result may be viewed as an \emph{input-output design separation principle} and plays a key role in characterizing the minimum number of feedback links required to ensure condition b) in Theorem \ref{Theorem2}.

\vspace{-0.3cm}

\begin{lemma}\label{Lemma5}
Let $\mathcal B(\mathcal X,\mathcal X,\mathcal E_{\mathcal X,\mathcal X})$ be the state bipartite graph. If $M^1$ and $M^2$ are two possible maximum matchings of $\mathcal B(\mathcal X,\mathcal X,\mathcal E_{\mathcal X,\mathcal X})$ with right-unmatched and left-unmatched vertices given by $(\mathcal V_R^1,\mathcal V_L^1)$ and $(\mathcal V_R^2,\mathcal V_L^2)$ respectively, then, exists a maximum matching $M^*$ of $\mathcal B(\mathcal X,\mathcal X,\mathcal E_{\mathcal X,\mathcal X})$ with right-unmatched and left-unmatched vertices given by  $(\mathcal V_R^1,\mathcal V_L^2)$. In particular, if $M^1$ has maximum top assignability index and $M^2$ has maximum  bottom assignability index, then  $M^*$ given by the above has both maximum top and bottom assignability index.\hfill $\diamond$
\end{lemma}
\vspace{-0.3cm}

Now referring to the discussion in Remark \ref{remarkStruct}, it follows that for a structurally controllable system, each non-top linked SCC in the state digraph must contain at least one state variable that is  reachable from an effective input, and, similarly,  by duality, we may conclude that for a structurally observable system, each non-bottom linked SCC in the state digraph must contain at least one state variable that  reaches an effective output. Now, noting that each SCC in the state digraph is either non-top linked or is reachable from a non-top linked SCC, we may further conclude that in a structurally controllable and observable configuration all state variables  are reachable from the effective inputs and, similarly, by duality,  all state variables reach effective outputs.

The above immediately yields the following property for the state-input-output digraph.

\vspace{-0.3cm}

\begin{proposition}
\label{prop:aux} 
Let $\bar B^e$ and $\bar C^e$ correspond to the collection of non-zero columns  of $\bar B$ and non-zero rows of $\bar C$ respectively, which correspond to the effective inputs and outputs  respectively.
If  $(\bar{A},\bar{B},\bar{C})$ is such that the pair $(\bar{A},\bar{B})$ is structurally controllable and $(\bar{A},\bar{C})$ is structurally observable, then,  the  effective inputs and  effective outputs constitute the non-top linked SCCs and non-bottom linked SCCs of $\mathcal D(\bar A,\bar B^e,\bar C^e)$, respectively. \hfill $\diamond$
\end{proposition}
\vspace{-0.3cm}

Now, recall Lemma \ref{Lemma2} and assume that a perfect matching of the state bipartite graph exists. Then there exists only one effective input and effective output (see Corollary \ref{minsolP1}, and more precisely  (13)), and  consequently  $\mathcal D(\bar A)$ is spanned by a disjoint union of cycles, which is sufficient to fulfill condition b) in Theorem \ref{Theorem2}. Hence, in that case,  by considering a single feedback link  from the effective output to the effective input  (which by Proposition \ref{prop:aux} are the unique non-top linked and non-bottom linked SCC in $\mathcal{D}(\bar{A},\bar{B}^{e},\bar{C}^{e})$ respectively), a single SCC is obtained with a feedback edge on it, additionally fulfiling condition a) in Theorem \ref{Theorem2}. In other words, an optimal information pattern  $\bar K$ has a single non-zero entry, more precisely, $\bar K_{ij}\neq 0$ corresponding to a feedback link between the $j$-th output (the only effective output) and the $i$-th input (the only effective input).

A more challenging case is encountered when a maximum matching is not perfect.  First, notice that by Remark \ref{rem:struc}, we need as many effective inputs/outputs as the number of right/left-unmatched vertices and, by Proposition \ref{prop:aux},  a feedback link is required for each effective input/output if condition a) of Theorem \ref{Theorem2} is to be satisfied. In other words, a lower bound for $\bar K$ in $\mathcal P_2$ may be obtained as
\vspace{-0.6cm}
\begin{equation}\label{lowerboundK}
\|\bar K\|_0\ge \max(m,1),
\vspace{-0.2cm}
\end{equation}
where $m$ is the number of right/left-unmatched vertices (that matches the number of effective inputs/outputs). This is so because, on one hand, using fewer effective inputs (respectively outputs)  will lead to loss of structural controllability (respectively observability), on the other hand, at least $m$ feedback links need to be considered to fulfil condition b) in Theorem \ref{Theorem2}.

We now show through a series of arguments that the lower bound in (12) is indeed achievable. First, we show that condition b) may be ensured by adding (appropriately) $\max(m,1)$ number of feedback links between effective outputs and inputs. To this end, in general, consider a minimal $(\bar B,\bar C)$ (see Corollary \ref{minsolP1}) to which a common matching exists (in the sense of Lemma \ref{Lemma5}), and by recalling Lemma \ref{lemmaB}, that such a maximum matching provides a decomposition of $\mathcal{D}(\bar{A})$ into a  disjoint union of cycles and state stems. Finally, we note that (see Corollary \ref{minsolP1}), distinct effective inputs (respectively outputs) are assigned to the roots (respectively tips) of such state stems rendering the latter to input-output stems. The above discussion is formalized as follows:

\vspace{-0.4cm}

\begin{proposition}\label{IOdecomp}
Let $\mathcal B$ be the state bipartite graph associated with $\mathcal D(\bar A)$, $\bar{B}$ and $\bar{C}$ be constructed as in Corollary \ref{minsolP1}. In addition, let $M^{\ast}$ be a (non-perfect) common maximum matching (in the sense of Lemma \ref{Lemma5}) such that its set of right-unmatched vertices $\mathcal{U}_{R}$ and left-unmatched vertices $\mathcal{U}_{L}$ correspond to the locations of effective inputs given by the collection of non-zero columns $\bar{B}^{e}$ of $\bar{B}$ and effective outputs given by the collection of non-zero rows $\bar{C}^{e}$ of $\bar{C}$ respectively.   Then, $\mathcal D(\bar A,\bar B^e,\bar C^e)$ is spanned by a disjoint union of cycles (composed exclusively of state vertices) and input-output stems. \hfill $\diamond$
\end{proposition}
\vspace{-0.3cm}

Specifically, note that the cycles provided by $M^{\ast}$ presents a covering of a subset of vertices in $\mathcal{D}(\bar{A})$, and, hence, as far as adding feedback edges are concerned to ensure condition b) in Theorem \ref{Theorem2}, only the state stems (which are input-output stems in $\mathcal{D}(\bar{A},\bar{B}^{e},\bar{C}^{e})$) provided by $M^{\ast}$ need to be considered (covered). Moreover, there are exactly $m$ such input-output stems and hence, in particular, by adding a feedback edge between the output and the input of each stem, we may obtain $m$ cycles that cover all the state vertices which belonged to state stems in the decomposition provided by $M^{\ast}$. As an immediate consequence, we note that $\max(m,1)$ feedback edges are sufficient in general to achieve condition b) in Theorem \ref{Theorem2}.

We remark that closing each input-output stem individually as explained above is not the only way to ensure condition b) in Theorem \ref{Theorem2} through $\max(m,1)$ feedback links. In particular, it may be possible to pair inputs and outputs belonging to different input-output stems using $\max(m,1)$ feedback edges and still satisfying the requirement b) in Theorem \ref{Theorem2}. Since, we are interested in obtaining all possible solutions to $\mathcal{P}_{2}$, we now provide a generic input-output pairing process,  which characterizes all possible input-output pairings with $\max(m,1)$ feedback edges that satisfy requirement b) in Theorem \ref{Theorem2}. Finally, we note that all such pairings may not satisfy condition a) in Theorem \ref{Theorem2}; however, we will show that there exists at least one such pairing which satisfies a) (and of course b)), thus establishing the achievability of the lower bound in (12). Each pairing satisfying both conditions of Theorem \ref{Theorem2}  will be referred to  as a \emph{mix-pairing}.

To this end, we introduce the \emph{IO-reachability bipartite graph}  $\mathcal B^{u,y}\equiv \mathcal B(\mathcal V_u,\mathcal V_y,\mathcal E_{u,y})$, with $\mathcal V_u$ and $\mathcal V_y$ denoting the set of effective inputs and effective outputs respectively,  and  an edge $(u,y) \in \mathcal E_{u,y}$, for $u\in \mathcal V_u$ and $y\in \mathcal V_y$ if and only if $u$ reaches $y$ in $\mathcal D(\bar A,\bar B,\bar C)$. Now, let $M^{u,y}$ be a maximum matching of $\mathcal{B}^{u,y}$, which is perfect by Proposition \ref{prop:aux}. Also,  observe that each edge in $M^{u,y}$ corresponds to a pair $(input, output)$, the latter being the root and tip of (possibly different) input-output stems (see Proposition \ref{IOdecomp}).

Formally, to characterize  all information patterns that satisfy Theorem \ref{Theorem2}-b),  we fix such a maximum matching $M^{u,y}$ and consider a partition of $M^{u,y}$ into disjoint  subsets  $\mathcal S_1,\cdots,\mathcal S_{|M^{u,y}|}$, i.e.,  $\bigcup_{i=1,\cdots,|M^{u,y}|} \mathcal S_i =M^{u,y}$ and $\mathcal S_i\cap \mathcal S_j = \emptyset$ for all $i\neq j$. Now, given $M^{u,y}$ and a partition, for each subset $\mathcal{S}_{i}$ we will denote by $\mathcal{F}(M^{u,y},\mathcal{S}_{i})$ the collection of all possible sets of $|\mathcal{S}_{i}|$ potential feedback edges whose incorporation augments the input-output stems involved in $\mathcal{S}_{i}$ into a single cycle. Specifically,
\vspace{-0.4cm}
\begin{align}
&\mathcal F(M^{u,y},\mathcal S_i)=\{\{f_1,\ldots,f_{|\mathcal S_i|}\}: \ e_1f_1\ldots e_{|\mathcal S_i|}f_{|\mathcal S_i|} \text{ describes a cycle in } \notag\\
&  \mathcal{B}(\mathcal{V}_{u},\mathcal{V}_{y},\mathcal{E}_{u,y}\cup\{f_{1},\cdots,f_{|\mathcal{S}_{i}|}\}) \text{ and }  e_1,\ldots,e_{|\mathcal S_i|} \text{ is an enumeration  of the edges in } \mathcal S_i \}. \label{defFcal}
\end{align}
\vspace{-1cm}

In particular, note  that $\mathcal F(M^{u,y},\mathcal S_i)$ has as many elements (sets of $|\mathcal S_i|$ edges) as the number of possible enumerations of the edges in $\mathcal S_i$.

In addition, we now define for each choice of subsets $\mathcal S_i$, $i=1,\cdots,|M^{u,y}|$, the set of sparsest information patterns (i.e., with  minimum $\|\bar K\|_0$) that satisfy Theorem \ref{Theorem2}-b) as
\vspace{-1.2cm}

\begin{align}
\mathcal K(M^{u,y},\mathcal S_1,\cdots,\mathcal S_{|M^{u,y}|})=\{\bar K\in \{0,1\}^{n \times n}: \ \bar{K}=\bar K (\mathcal E_{\mathcal Y,\mathcal U})& \text{ where } \mathcal E_{\mathcal Y,\mathcal U} =\bigcup_{i=1,\ldots,|M^{u,y}|}\mathcal E_i\ ,\notag \\
& \text{ with } \mathcal E_i \in \mathcal F(M^{u,y},\mathcal S_i)\},\label{defKcal}
\end{align}
where $\bar K (\mathcal E_{\mathcal Y,\mathcal U})$  denotes the structural matrix $\bar K$ whose non-zero entries correspond to the edge set $\mathcal{E}_{\mathcal{Y},\mathcal{U}}$ in the digraph representation $\mathcal{D}(\bar{A},\bar{B},\bar{K},\bar{C})$ (see definition of $\mathcal D(\bar A,\bar B,\bar K,\bar C)$).

Notice that not all information patterns in $\mathcal K(M^{u,y},\mathcal S_1,\cdots,\mathcal S_{|M^{u,y}|})$ satisfy both the conditions in Theorem \ref{Theorem2}, see Figure~\ref{fig:mixpairing}. However, the following result holds.

\vspace{-0.2cm}

\begin{proposition}
\label{prop:aux2}
Let $\mathcal D (\bar A,\bar B,\bar C)$  be the systems digraph and $\mathcal B^{u,y}$ its IO-reachability bipartite graph.  Then, the following conditions hold:


\begin{itemize}
\item[(1)]  If $\bar K \in \mathcal K (M^{u,y},\mathcal S_1,\ldots,\mathcal S_{|M^{u,y}|})$ for some choice of maximum matching $M^{u,y}$ of $\mathcal{B}^{u,y}$ and subsets $\mathcal{S}_{i}\subset M^{u,y}$, $i=1,\cdots,|M^{u,y}|$,   with $\mathcal S_i\cap \mathcal S_j=\emptyset$ for $i\neq j$ and $\bigcup_{i=1,\ldots,|M^{u,y}|}\mathcal S_i = M^{u,y}$, then, $\bar K$ is such that $\mathcal D(\bar A,\bar B,\bar K,\bar C)$ satisfies condition b) in Theorem \ref{Theorem2}.
\item[(2)] For a maximum matching $M^{u,y}$ and  the specific choice of $S_i = M^{uy}$ for a particular $i$ in $\{1,\cdots,|M^{u,y}|\}$, i.e., $\mathcal{S}_{j}=\emptyset$ for $j\neq i$,  any  information pattern $\bar{K}\in\mathcal K(M^{u,y},\mathcal S_1,\cdots,\mathcal S_{|M^{u,y}|})$ is such that $D(\bar A,\bar B,\bar{K}, \bar{C})$ satisfies both the conditions in Theorem \ref{Theorem2};
\item[(3)] If $\bar{K}$ is a sparsest  information pattern (i.e., with minimal $\|\bar{K}\|_{0}$) such that $\mathcal D(\bar A,\bar B,\bar K,\bar C)$ satisfies both the conditions in Theorem \ref{Theorem2}, then there exists a maximum matching $M^{u,y}$ of $\mathcal B^{u,y}$ and subsets $\mathcal S_1,\cdots,\mathcal S_{|M^{u,y}|}\subset M^{u,y}$, with $\mathcal S_i\cap \mathcal S_j=\emptyset$ for $i\neq j$ and $\bigcup\limits_{i=1,\ldots,|M^{u,y}|}\mathcal S_i = M^{u,y}$, such that $\bar{K}\in\mathcal K(M^{u,y},\mathcal S_1,\cdots,\mathcal S_{|M^{u,y}|})$.
\hfill $\diamond$
\end{itemize}
\end{proposition}
\vspace{-0.2cm}



%

In particular, it follows by Proposition \ref{prop:aux2}-(2)  that there exists $\bar K$ such that $\|\bar K\|_0=\max(m,1)$, which together with $(\bar B,\bar C)$ (constructed as in Corollary~\ref{minsolP1}) is a solution to $\mathcal P_2$.

Additionally, we denote by  $\mathcal K_{mix}(M^{u,y},\mathcal S_1,\cdots,\mathcal S_{|M^{u,y}|})$   the subset of sparsest information patterns in  $\mathcal K(M^{u,y},\mathcal S_1,\cdots,\mathcal S_{|M^{u,y}|})$ that satisfy Theorem \ref{Theorem2}-a)  for a given maximum matching $M^{u,y}$ of $\mathcal B^{u,y}$ and subsets $\mathcal S_1,\cdots,\mathcal S_{|M^{u,y}|}$ of $M^{uy}$. Therefore, all the sparsest information patterns that satisfy problem $\mathcal P_2$ upon a given construction of $(\bar B,\bar C)$ as in Corollary~\ref{minsolP1}, are   characterized by 
\vspace{-0.5cm}
\begin{equation}
\label{rem:mix_pairing}
\mathcal K_{mix} = \bigcup_{M^{u,y}; \ \mathcal S_1,\cdots,\mathcal S_{|M^{u,y}|}\subset M^{u,y}}\mathcal K_{mix}(M^{u,y},\mathcal S_1,\cdots,\mathcal S_{|M^{u,y}|}).
\end{equation}

For notational convenience, we refer to any $\bar K \in \mathcal K_{mix}$ as a \emph{mix-pairing} of $\mathcal D(\bar A,\bar B,\bar C)$. Figure~\ref{fig:mixpairing} depicts two examples of mix-pairing.

\begin{remark}
\label{rem:mix_pairing} 
Computing  a mix-pairing instance incurs in polynomial complexity  because  it reduces to finding a maximum matching $M^{uy}$ and subsequently resorting to the construction proposed in Proposition~\ref{prop:aux2}-(2).\hfill $\diamond$
\end{remark}

\vspace{-0.5cm}

\begin{figure}[htb]
\centering
\includegraphics[scale=0.33]{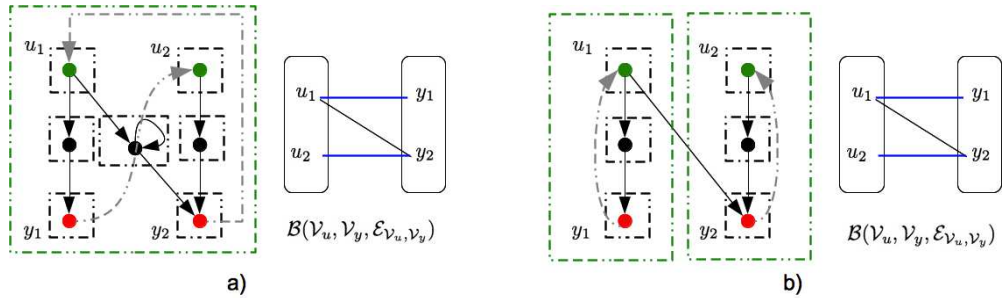}
\caption{a) A mix-pairing $\bar K \in \mathcal F(M^{u,y},\mathcal S_1,\mathcal S_2)$
with $\mathcal S_1=M^{u,y}$ and $\mathcal S_2=\emptyset$ is depicted by gray dashed edges, all the SCCs in the original digraph are  depicted by black dashed boxes, whereas the SCCs associated with the new digraph $\mathcal D(\bar A,\bar B,\bar K,\bar C)$ are depicted by green dashed boxes. The blue edges  represent a possible maximum matching associated with the IO-reachability bipartite graph.  b) A mix-pairing $\bar K  \in \mathcal F(M^{u,y},\mathcal S_1,\mathcal S_2)$ with $\mathcal S_1,\mathcal S_2$ each comprising a single edge from $M^{u,y}$ is depicted by gray dashed edges. Remark that the mix-pairing used in b) cannot be used in a), otherwise, the digraph $\mathcal D(\bar A,\bar B,\bar K,\bar C)$ will comprise of three SCCs, with one SCC (corresponding to the node with the self-loop) having no feedback edge in it, thus violating Theorem \ref{Theorem2}-a).}
\label{fig:mixpairing}
\end{figure}

\vspace{-.5cm}

Finally, we state the main result of this section.

\vspace{-0.4cm}

\begin{theorem}[Solution to $\mathcal P_2$]\label{Theorem8}
 The triple  $(\bar B,\bar K,\bar C)$ is a solution to $\mathcal P_2$ \underline{if and only if} the pair $(\bar B,\bar C)$ is a minimal solution to $\mathcal P_1$ (constructed as in Corollary \ref{minsolP1}) and $\bar K$  corresponds to   a mix-pairing of $\mathcal D(\bar A,\bar B,\bar C)$. \hfill $\diamond$
\end{theorem}
\vspace{-0.6cm}

\begin{remark}\label{ShreyasWork}
We contrast our CC selection results to that of \cite{DBLP:journals/jsac/PajicMPS13}.   In \cite{DBLP:journals/jsac/PajicMPS13}, for discrete time linear invariant systems, given a set of effective inputs and outputs, a pairing process similar to the construction in Proposition \ref{prop:aux2}-(2)  is provided. In contrast, we provide all possible sparsest information patterns for the CC selection problem in the continuous time scenario. Moreover, from a technical standpoint, due to the discrete time treatment, the construction in \cite{DBLP:journals/jsac/PajicMPS13} needs to ensure that condition a) in Theorem \ref{Theorem2} is satisfied only (the uncontrollable/unobservable modes at zero pose no concern in the current context for the discrete time setting), whereas, the requirement to satisfy both conditions a) and b) simultaneously in the continuous time setting adds a layer of technical complexity in our construction and analysis. \hfill $\diamond$
\end{remark}
}


\vspace{-0.9cm}

\section{Algorithmic Procedure and Complexity Analysis}\label{complexAnalysis}
\vspace{-0.2cm}

We now provide an efficient algorithmic procedure to compute a minimal feasible dedicated input configuration $\mathcal S_u$,  described in Algorithm 1. Briefly, Algorithm 1 consists in finding a maximum matching with maximum top assignability and its associated set of right-unmatched vertices. Then, by Theorem \ref{MFDIC}, a minimal feasible dedicated input configuration may be obtained by assigning (dedicated) inputs to state variables corresponding to the right-unmatched vertices of the maximum matching and an additional set of $\beta-\alpha$  state variables each of which belongs  to a distinct non-assigned non-top linked SCC. 
\vspace{-0.5cm}

\begin{algorithm}[th]
\small
\KwIn{$\mathcal D(\bar A)=(\mathcal X,\mathcal E_{\mathcal X,\mathcal X})$}
\KwOut{Minimal feasible dedicated input configuration $\mathcal S_u$}

$\mathbf{Step  \ 1.}$ Determine the  non-top linked SCCs  $\mathcal{N}^T_i,\ i\in \mathcal I\equiv\{1,\cdots, \beta\}$, of $\mathcal D(\bar A)$ and denote its collection by $\mathcal N$.

$\mathbf{Step \ 2.}$ Consider a weighted bipartite graph $\mathcal B(\mathcal X,\mathcal X\cup\mathcal {I},\mathcal E_{\mathcal X,\mathcal X}\cup \mathcal E_{\mathcal I,\mathcal X})$ where $\mathcal E_{\mathcal I,\mathcal X}=\{(i,x_j): \ x_j\in \mathcal N^T_i\}$ and each variable $i\in\mathcal I$ is a slack variable with an edge to each variable in the non-top linked SCC $\mathcal N^T_i$. The associated cost is given as follows: each edge that does not belong to the bipartite graph has infinite cost, each edge in $\mathcal E_{\mathcal X,\mathcal X}$ has unitary cost and each edge in $\mathcal E_{\mathcal I,\mathcal X}$ has  a  cost  of two.

$\mathbf{Step \ 3.}$ Let $M'$ be the maximum matching incurring in the minimum cost of the weighted bipartite graph presented in step 2 and $\mathcal U_R$ be the corresponding set of right-unmatched vertices (which may comprise varibles in $\mathcal I$). Then, $\mathcal U_R^*=(\mathcal U_R\cap \mathcal X) \cup \Theta$, where $\Theta = \{ x_j: (i,x_j) \in M', i\in \mathcal I\}$, is the set of right-unmatched vertices associated with a maximum matching of the state bipartite graph $\mathcal B(\mathcal X,\mathcal X,\mathcal E_{\mathcal X,\mathcal X})$ with maximum top assignability.

$\mathbf{Step \ 4.}$ Set $\mathcal S_u=\mathcal U^*_R\cup \mathcal A_u^c$ where $\mathcal A_u^c$ consists of a union of state variables formed by selecting a single state variable from each non-assigned SCC $\mathcal N^T_j$ (i.e., $(\mathcal N^T_j, \cdot)\notin M^*$).

\caption{Computing a minimal feasible dedicated input configuration $\mathcal S_u$}
\label{mainAlgorithmNew}
\end{algorithm}

\vspace{-0.5cm}
The  next set of results establishes  the correctness and analyzes the implementation complexity  of Algorithm 1. 

\vspace{-0.3cm}

\begin{theorem}[Correctness and Computational Complexity of Algorithm \ref{mainAlgorithmNew}]\label{TheoremCorrectness}
Algorithm \ref{mainAlgorithmNew} is correct, i.e., its execution provides a minimal feasible dedicated input configuration $\mathcal S_u$. Furthermore, it generates a minimal feasible dedicated input configuration $\mathcal S_u$, with complexity  $\mathcal O(|\mathcal X|^{3})$, where $|\mathcal X|$ denotes the number of state vertices in $\mathcal D(\bar A)$. In addition,    the minimum number of dedicated inputs to ensure structural controllability of the system is given by  $|\mathcal S_u|$ and its computation incurs in the same complexity.
\hfill $\diamond$
\end{theorem}
\vspace{-0.3cm}

Given that by formulation the problem of computing a minimum feasible dedicated input configuration is a combinatorial optimization problem, the polynomial complexity construction provided in Algorithm 1 is especially helpful in the context of large-scale systems. To emphasize further, even assuming that the minimum number $p$ of dedicated inputs required $p$ is known, a naive combinatorial search over  $\binom{|\mathcal{X}|}{p} $ possible configuration choices (and verifying if each of them is feasible or not, which may be achieved using an algorithm of quadratic complexity in the number of state variables  \cite{dionSurveyKyb}) may not be feasible in large-scale scenarios; in fact, if $p$ grows with $|\mathcal{X}|$, such a combinatorial search procedure may incur exponential complexity.

Noting that solutions to $\mathcal{P}_{1}$ and $\mathcal{P}_{2}$ may be obtained by simplistic constructions once a minimal feasible dedicated input configuration is provided (see Theorem \ref{solP1} and Theorem \ref{Theorem8}).
\vspace{-0.3cm}

\begin{corollary}
There exist  $\mathcal O(|\mathcal X|^{3})$ complexity procedures for computing solutions to $\mathcal{P}_{1}^{d}$, $\mathcal{P}_{1}$ and $\mathcal{P}_{2}$,  where $|\mathcal X|$ denote the number of state vertices in $\mathcal D(\bar A)$.
\hfill $\diamond$
\label{Corollary5}
\end{corollary}

\vspace{-0.9cm}

\section{Illustrative Example}\label{illustrativeexample}
\vspace{-0.2cm}

The following example illustrates the procedure to obtain a solution to  $\mathcal P_1$ and $\mathcal P_2$ following  the sequence of results presented in Sections III-V. First we compute a solution to (3a) and (3b) to be used in computing a solution to $\mathcal{P}_{1}$. Then,  we use a solution of $\mathcal P_1$ to compute a solution to $\mathcal P_2$.

\begin{figure}[htb]
\centering
\includegraphics[scale=0.26]{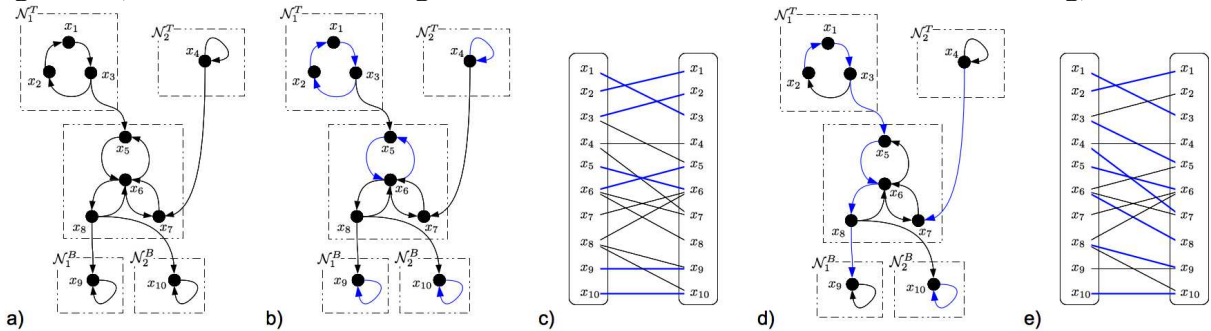}
\caption{a) The original digraph $\mathcal{D}(\bar A)=(\mathcal X,\mathcal E_{\mathcal X,\mathcal X})$ where  the dashed boxes denote the SCCs of $\mathcal D(\bar A)$, the black vertices the state vertices, and the arrows the directed edges.; b) The digraph $\mathcal{D}(\bar A)=(\mathcal X,\mathcal E_{\mathcal X,\mathcal X})$ where  the edges in blue correspond to a maximum matching, associated with $\mathcal B(\mathcal X,\mathcal X,\mathcal E_{\mathcal X,\mathcal X})$, depicted in c); d)  The  digraph $\mathcal{D}(\bar A)=(\mathcal X,\mathcal E_{\mathcal X,\mathcal X})$ where the edges in blue belong to a maximum matching with maximum top and bottom assignability, associated with $\mathcal B(\mathcal X,\mathcal X,\mathcal E_{\mathcal X,\mathcal X})$, depicted in e). }
\label{fig:ex1}
\end{figure}

Consider the directed graph  $\mathcal{D}(\bar A)=(\mathcal X,\mathcal{E}_{\mathcal X,\mathcal X})$  depicted in Figure \ref{fig:ex1}-a) and compute a maximum matching associated with $\mathcal{B}(\mathcal{X},\mathcal{X},\mathcal{E}_{\mathcal{X},\mathcal{X}})$.
Figure \ref{fig:ex1}-b) represents in blue, the edges belonging to  the maximum matching $M^*=\{(x_2,x_1),$ $(x_1,x_3),(x_3,x_2), (x_4,x_4),(x_5,x_6),(x_6,x_5),$ $(x_9,x_9),(x_{10},x_{10})\}$ (note that,  in general, the maximum matching is not unique, for instance,  Figure \ref{fig:ex1}-c) shows in blue the edges associated with a different maximum matching).

\vspace{-0.7cm}

\vspace{-0.9cm}

\underline{Solution to (3a)} From $M^*$, the set of right-unmatched vertices is $\mathcal V=\{x_8,x_7\}$ and therefore by  Theorem~\ref{Theorem4} we have $m=2$. Moreover,  because $\mathcal D(\bar A)$ has two non-top linked SCCs ($\beta$=2 in  Theorem~\ref{Theorem4}). To find $\alpha$ in Theorem~\ref{Theorem4}, which is defined to be the maximum top assignability index of $\mathcal B(\mathcal X,\mathcal X,\mathcal E_{\mathcal X,\mathcal X})$, we have to check, in general, if there is another set of right-unmatched vertices (corresponding to another maximum matching of $\mathcal B(\mathcal X,\mathcal X,\mathcal E_{\mathcal X,\mathcal X})$) that covers (i.e., has elements in) a larger number of non-top linked SCCs\footnote{Note that the purpose of this example is to illustrate the technical results established in the paper. To that end and since the example system is small, the various constructions are mostly carried out by hand. For larger systems, the systematic algorithmic procedures given in Section VI of the paper can be used.}. This can be done as follows: check if there is a vertex in the  non-top linked SCC  $N^T_1$ or $N^T_2$ that could replace $x_8$ as a unmatched vertex, i.e.,  verify if one of the following sets $\{x_1,x_7\},\{x_2,x_7\},\{x_3,x_7\},\{x_4,x_7\}$ corresponds to the set of right-unmatched vertices of a maximum matching (which can be implemented using Proposition~\ref{prop:aux}  in Section II). The set $\{x_1,x_7\}$ does not correspond to a set of right-unmatched vertices, i.e., considering the state bipartite graph without the edges ending in $x_1,x_7$, by computing its maximum matching, we obtain always an associated set of right-unmatched vertices that strictly contains  $\{x_1,x_7\}$, for instance,  $\{x_1,x_7,x_8\}$. Now, by considering $\{x_2,x_7\}$ and proceeding similarly, we conclude that it corresponds to a possible set of right-unmatched vertices where the non-top linked SCC $\mathcal N^T_1$ is assigned.
Therefore, we need to determine if $\mathcal N^T_2$ is also assignable, which reduces to checking if the set $\{x_2,x_4\}$ corresponds to the set of right-unmatched vertices associated with a maximum matching of the bipartite graph obtained by removing the edges ending in $x_{2}$ and $x_{4}$ from $\mathcal{B}(\mathcal{X},\mathcal X,\mathcal E_{\mathcal X,\mathcal X})$.  Indeed, this is the case and hence both $\mathcal{N}^T_1$, $\mathcal N^T_2$ are assignable, i.e.,  $\alpha=2$ in Theorem~\ref{Theorem4}. The edges corresponding to a maximum matching with maximum top  assignability index are depicted in  Figure \ref{fig:ex1}-c). Hence, we need $p=m+\beta-\alpha=2+2-2=2$ dedicated inputs. To systematically compute $\alpha$ one can compute the maximum matching of the set computed in Algorithm 1.

\begin{figure}[htb]
\centering
\includegraphics[scale=0.25]{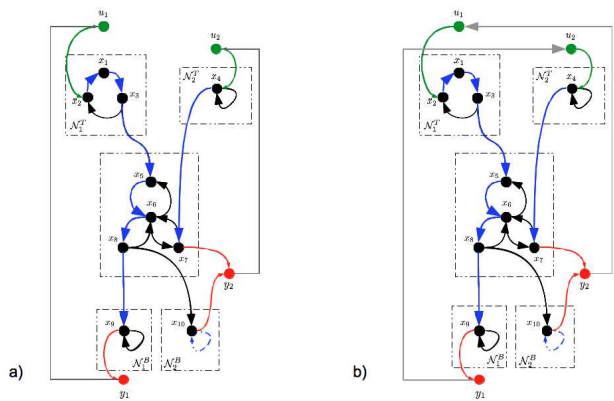}
\caption{In a) and b) we depict two different solutions to $\mathcal P_2$. In blue we can identify the directed edges associated with matched edges in a maximum matching, in green the input edges and in red the output edges. From the solid colored lines we can easily identify the cycles comprising input-output stems: in a) the feedback edge in gray closes the loop with an input-output stem whereas in b) there is only one cycle comprising two input-output stems.}
\label{fig:ex2}
\end{figure}

Therefore, a particular minimal  feasible input configuration (see Theorem~\ref{MFDIC}) is $\{x_2,x_4\}$. Remark that this can also be determined using Algorithm 1. In fact, $\{x_2,x_4\}$ is the only possible minimal feasible input configuration. Hence, by assigning an input to each of $\{x_{2},x_{4}\}$, we obtain a structurally controllable system.
This choice is depicted in Figure \ref{fig:ex2}-a), but notice that the choice of the index of the effective inputs ($u_1,u_2$) is arbitrary, corresponding to a permutation of the columns of the matrix $\bar B$.  The structure of $\bar B$ is given by $\bar B_{2,1}=1$, $\bar B_{4,2}=1$ and zero elsewhere. Moreover, (3a) has  a total of $\|\bar B\|_0=2$ edges between inputs and states  (depicted by the green edges in Figure \ref{fig:ex2}-a)).

\underline{Solution to (2a)}: In  this case the solution of (3a) is also a solution to (2a) under the additional constraint of minimum number of effective inputs.  By  Theorem~\ref{sol2a} (and also Corollary~\ref{minsol2a}),  the structure of $\bar B$ is given by $\bar B_{2,1}=1$, $\bar B_{4,2}=1$ and zero elsewhere. Moreover, (2a) has  a total of $\|\bar B\|_0=2$ edges between inputs and states (depicted by the green edges in Figure \ref{fig:ex2}-a)). 

 \underline{Solution to (3b)}: To compute the minimum number of state variables to be assigned with outputs to have structural observability, we proceed with the same reasoning as above but using $\mathcal D (\bar A^T)$ as the system digraph, which corresponds to reversing the directions of the edges of the original digraph $\mathcal D(\bar A)$. Thus, $m'=2$, $\beta'=2$ and $\alpha'=1$, where, in terms of the original digraph $\mathcal D (\bar A)$, $m'$ denotes the number of left-unmatched vertices of a maximum matching associated with the bipartite graph $\mathcal B(\mathcal X,\mathcal X,\mathcal E_{\mathcal X,\mathcal X})$, $\beta'$ denotes the number of non-bottom linked SCCs and $\alpha'$ the maximum bottom assignability index. The edges corresponding to a maximum matching with maximum  bottom assignability index are depicted in  Figure \ref{fig:ex1}-c). Hence, a total of $p'=m'+\beta'-\alpha'=2+2-1=3$ state variables to ensure structural observability are required.
The minimal feasible output configuration is  given by $\{x_9,x_{10},x_7\}$, by the dual of Theorem~\ref{MFDIC} applied to  structural ouput design.
\vspace{-0.7cm}

\vspace{-0.7cm}

 \underline{Solution to (2b)} : By the dual of Theorem~\ref{sol2a}, the  solution of (2b) consists of two outputs (say $y_1,y_2$) that are required to measure the state variables associated with the left-unmatched vertices, i.e., $\{x_7,x_9\}$. In addition, a new output or one of the previously assigned should be assigned to $x_{10}$. For the case that $y_1$ or $y_2$ is connected to $x_{10}$ we have the minimal output solution (given by the dual of Corollary~\ref{minsol2a} applied to  structural design), depicted in Figure \ref{fig:ex2}-a). Therefore, the structure of $\bar C$ associated with the previous choice is:  $\bar C_{1,9}=1$, $\bar C_{2,10}=1$, $\bar C_{2,7}=1$ and zero elsewhere. Hence, (2b) incurs in the minimum of $\|\bar C\|_0=3$ (depicted by the red edges in Figure \ref{fig:ex2}-a)).

\underline{Solution to $\mathcal P_1$} Merge the solutions of (2a) and (2b) and the result follows by Theorem~\ref{solP1}. In particular, we also have the minimal solution to $\mathcal P_1$, as described in Corollary~\ref{minsolP1}.


\underline{Solution to $\mathcal P_2$} Invoking Theorem~\ref{Theorem8} we can construct a solution to $\mathcal P_2$. Consider $\bar B$ and $\bar C$ previously constructed to solve $\mathcal P_1$, which are also minimal solutions (in the sense of  Corollary~\ref{minsolP1}). First, notice that a common maximum matching exists (as described by Lemma~\ref{Lemma5}), which is depicted by the blue edges in Figure \ref{fig:ex2}. In fact, notice that $u_1\rightarrow x_2\rightarrow x_1 \rightarrow x_3 \rightarrow x_5 \rightarrow x_6 \rightarrow x_8 \rightarrow x_9 \rightarrow y_1$ defines an input-output stem, as well as $u_2 \rightarrow x_4 \rightarrow x_7 \rightarrow y_2$. Together with the feedback edges depicted by gray edges in Figure \ref{fig:ex2}, the input-output stems are covered by cycles, as specified in (13)-(14). Nevertheless, notice that although both Figure \ref{fig:ex2}-a) and Figure \ref{fig:ex2}-b) satisfy condition b) in Theorem 2, Figure \ref{fig:ex2}-a) does not fulfill condition a) of Theorem 2 since $x_{10}$ does not belong to an SCC with a feedback edge on it. Hence, the only feasible solution, i.e., the only mix-pairing is given by  $\bar K$ such that  $\bar K_{1,2}=\bar K_{2,1}=1$ and zero elsewhere, depicted in Figure \ref{fig:ex2}-b).

\vspace{-0.6cm}

\section{Conclusions and Further Research}\label{secConclusion}
\vspace{-0.2cm}

In this paper we have proposed several novel graph-theoretic characterizations of all possible solutions to the following structural system design problems: (1) the I/O selection problem, which studies the sparsest configurations of state variables that are to be actuated/measured by inputs/outputs (either dedicated or non-dedicated) to achieve structural controllability/observability, and (2) the CC problem, which studies the sparsest information patterns or equivalently minimum number and configurations of feedback edges required between inputs and outputs such that the closed loop system is free of structurally fixed modes. Our framework is unified, in that, it presents a joint solution to the above problems. Furthermore, we have provided algorithms of polynomial complexity (in the number of state variables) for systematically computing solutions to these problems.  A natural  direction for future research consists of extending the current framework to cope with I/O and CC selection with heterogeneous costs, i.e., in which different state variables incur different actuation/measurement costs and the communication cost (or cost of incorporating a feedback edge) varies from one input-output pair to another.

\vspace{-0.6cm}


%
\section*{APPENDIX A}

{
\setstretch{1}
\small

In the sequel, given two digraphs $\mathcal D^1$ and $\mathcal D^2$,  $\mathcal V(\mathcal D^i)$ and $\mathcal E(\mathcal V(\mathcal D^i))$, for $i=1,2$,  denote the set of vertices  and the set of  edges with both linked vertices in $\mathcal V(\mathcal D^i)$ respectively. 
In addition, the notation $\mathcal{D}^{1}\setminus\mathcal{D}^{2}$ denotes the digraph $\mathcal D^1\backslash \mathcal D^2=(\mathcal V^{\mathcal D^1\backslash \mathcal D^2}\equiv \mathcal V(\mathcal D^1)\backslash \mathcal V (\mathcal D^2), \mathcal E(V^{\mathcal D^1\backslash \mathcal D^2}))$.

\textit{Proof of Lemma \ref{lemmaA}} : 

Let $\mathcal D'=(\mathcal V',\mathcal E')$ be a digraph where each vertex $v^i\in \mathcal V'$ with $i=1,\ldots,c_{\mathcal D}$ corresponds to a cycle $\mathcal C^i$ in the original digraph $\mathcal D$, and an edge $(v^{i'},v^{i''})\in\mathcal E'$ exists  between two vertices  $v^{i'},v^{i''}\in\mathcal V'$ if and only if there exists an edge (in the original digraph $\mathcal{D}$) from a state vertex in $\mathcal C^{i'}$ to a state vertex in $\mathcal C^{i''}$. Now, notice that $\mathcal D'$ is an SCC, since the original $\mathcal D$ is an SCC. In particular, it follows that, for each vertex $v^{\prime}\in\mathcal{V}^{\prime}$, there exists a DAG $\mathcal{D}^{\prime\prime}$, with $v^{\prime}$ as its root, that spans $\mathcal{D}^{\prime}$. Thus,  the claim follows by noting that such a spanning DAG  in $\mathcal{D}^{\prime}$ corresponds to a chain in $\mathcal{D}$. \hfill $\blacksquare$

%
%

\textit{Proof of Lemma \ref{lemmaB}} : Let $\mathcal{C}^{\prime}_{i}$, $i=1,\cdots,\gamma$, be any subcollection of cycles such that $\mathcal{C}^{\prime}_{i}\in\mathcal{N}_{i}^{T}$ for all $i$. By applying Lemma \ref{lemmaA} to each SCC of $\mathcal{D}$, note that there exists a collection of (disjoint) chains each of which spans a distinct SCC of $\mathcal{D}$. Moreover, by Lemma \ref{lemmaA}, the spanning chain in the $i$-th non-top linked SCC, $i=1,\cdots,\gamma$, 
 may be constructed such that its first element is the cycle $\mathcal{C}^{\prime}_{i}$. Finally, note that each chain spanning an SCC that is non-top linked has an incoming edge to at least one of its vertices from (at least) one of the $\gamma$ chains spanning the non-top linked SCCs. In particular, by merging each chain in a top-linked SCC with exactly one chain in a non-top linked SCC, we obtain a new collection of disjoint $\gamma$ chains which span $\mathcal{D}$ and satisfy the requirements of Lemma \ref{lemmaB}.
 \hfill $\blacksquare$

%

\textit{Proof of Lemma \ref{Lemma2}} :
First observe from the definition, the subgraph $\mathcal D^*$ spans $\mathcal D(\bar A)$.

We next  prove the minimality of the decomposition (through state-stems and cycles) achieved by $M^{\ast}$. To this end, note that the following generic digraph properties may be verified from the definitions:

(1)  For $\mathcal{C}\subset\mathcal{E}_{\mathcal{X},\mathcal{X}}$, the digraph $\mathcal D^{\mathcal C}\equiv (\mathcal X,\mathcal{C})$ corresponds to a spanning decomposition of $\mathcal{D}(\bar{A})$ into disjoint subgraphs  of state stems and cycles \emph{if and only if } the edges in $\mathcal{C}$ define a matching $M^{\mathcal{C}}$ (not necessarily maximum) for $\mathcal{B}(\mathcal{X},\mathcal{X},\mathcal{E}_{\mathcal{X},\mathcal{X}})$.

(2) The root of a state stem belonging to a matching $M$ of $\mathcal{B}(\mathcal{X},\mathcal{X},\mathcal{E}_{\mathcal{X},\mathcal{X}})$ is necessarily a right-unmatched vertex w.r.t. $M$.

A direct consequence of the above properties is that the edges in $M^{\ast}$ constitute a spanning decomposition of $\mathcal{D}(\bar{A})$ into disjoint state stems and cycles. We  now establish the desired minimality of $M^{\ast}$ by contradiction. Assume, on the contrary, there exists a spanning decomposition $\mathcal D^{\mathcal{C}}$ of $\mathcal{D}(\bar{A})$ into disjoint subgraphs  of state stems and cycles, such that, $\mathcal{C}$ consists of strictly fewer state stems than $M^{\ast}$. Then, by property (1) above, the corresponding matching $M^{\mathcal{C}}$ consists of fewer state stems than $M^{\ast}$. Since, by property (2) above, each state-stem in $M^{\mathcal{C}}$ corresponds to a unique right-unmatched vertex (and similarly for $M^{\ast}$), we conclude that $M^{\mathcal{C}}$ consists of strictly fewer right-unmatched vertices than $M^{\ast}$. This clearly contradicts the fact, that, $M^{\ast}$ is a maximum matching. Hence, the desired assertion follows.
\hfill $\blacksquare$

%
%
%
%
%

\textit{Proof of Theorem \ref{FDIC}} :
[$\Rightarrow$] Let $\mathcal{S}_{u}$ be a feasible dedicated input configuration. Then, by Theorem 1, there exists a decomposition of $\mathcal{D}(\bar{A})$ into state cacti, where the root of each cactus is a vertex in $\mathcal{S}_{u}$. Now, further decomposing the state cacti collection into disjoint state stems and cycles (which necessarily span $\mathcal{D}(\bar{A})$), it follows that the roots of the resulting state stems are vertices in $\mathcal{S}_{u}$. Then, by property (1) in the proof of Lemma \ref{Lemma2}, we may conclude that  $\mathcal S_u$ corresponds to the set of right-unmatched vertices of a matching $M$ (not necessarily maximum) of $\mathcal B$.  Finally, note that there exists a maximum matching $M^{\ast}$, such that the set $\mathcal{U}_{R}$ of its right-unmatched vertices is a subset of the set of right-unmatched vertices of $M$. (Note that any matching can be augmented to a maximum matching, see~\cite{Cormen:2001:IA:580470}). It follows immediately that $\mathcal{U}_{R}\subset\mathcal{S}_{u}$.

Now assume on the contrary that there exists no $\mathcal A_u\subset \mathcal S_u$ containing at least one state vertex from each non-top linked  SCC. Then, there is a non-top linked SCC without any of its vertices in $\mathcal S_u$, which implies, by definition, that none of its vertices are  reachable from any state vertex in $\mathcal S_u$. Hence, in particular, such a non-top linked SCC cannot belong to any state cacti decomposition of $\mathcal{D}(\bar{A})$ with roots in $\mathcal{S}_{u}$, which contradicts that $\mathcal{S}_{u}$ is a feasible dedicated input configuration.

[$\Leftarrow$] We show that $\mathcal D(\bar A)$ can be decomposed into two partitions $\mathcal D_{\mathcal U_R}$ and $\mathcal D_{\mathcal A_u}$  that are spanned by disjoint union of state cacti with  roots in $\mathcal U_R$ and $\mathcal A_u$ respectively. Let $\mathcal D_{\mathcal U_R}=(\mathcal V_{\mathcal U_R},\mathcal E_{\mathcal U_R})$ where $\mathcal V_{\mathcal U_R}$ is composed of all state vertices that are reachable in $\mathcal D(\bar A)$ from the state vertices in $\mathcal U_R$,  and $\mathcal E_{\mathcal U_R}$ comprises the edges of $\mathcal D(\bar A)$ with both endpoints in $\mathcal V_{\mathcal U_R}$, and let  $\mathcal D_{\mathcal A_u}=\mathcal D(\bar A)\backslash \mathcal D_{\mathcal U_R}$. We now show that each partition is spanned by state cacti (with roots in $\mathcal U_R$ and $\mathcal A_u$ respectively) which would imply that $\mathcal D(\bar A)$ is  spanned by state cacti with roots in $\mathcal{S}_{u}$.

\indent $\bullet$ $ \mathcal D_{\mathcal U_R}$ is spanned by state cacti with roots in $\mathcal U_R$: By Lemma \ref{Lemma2}, it follows that $\mathcal D(\bar A)$ is spanned by a disjoint union of state stems  $\mathcal S$ with roots in $\mathcal U_R$ and a disjoint union of cycles $\mathcal C$. In particular, by construction, $\mathcal D_{\mathcal U_R}$ is spanned by $\mathcal S$ and a subset of the cycles that span $\mathcal D(\bar A)$, denoted by $\mathcal C'\subset \mathcal C$. The digraph $\mathcal D_{\mathcal U_R}\backslash \mathcal S$ is then  spanned by $\mathcal C'$ (a disjoint union of cycles) and hence, by Lemma 2, it follows that $\mathcal D_{\mathcal U_R}\backslash \mathcal S$ is spanned by $\gamma$ disjoint chains, where $\gamma$ denotes the number of non-top linked SCCs in $\mathcal{D}_{\mathcal{U}_{R}}\setminus\mathcal{S}$; moreover, each such chain may be constructed to have as its first element a cycle belonging to a distinct non-top linked SCC of $\mathcal{D}_{\mathcal{U}_{R}}\setminus\mathcal{S}$. 

Finally, note that, by construction of $\mathcal{D}_{\mathcal{U}_{R}}$, each non-top linked SCC in  $\mathcal{D}_{\mathcal{U}_{R}}\setminus\mathcal{S}$ has at least one incoming edge from the set $\mathcal{S}$. This implies, in particular, that the state stems $\mathcal{S}$ in $\mathcal{D}_{\mathcal{U}_{R}}$ connect to the chains spanning $\mathcal{D}_{\mathcal{U}_{R}}\setminus\mathcal{S}$ and hence, $\mathcal{D}_{\mathcal{U}_{R}}$ is  spanned by a disjoint union of state cacti with its roots in $\mathcal U_R$. 

\indent $\bullet$ $ \mathcal D_{\mathcal A_u}$ is spanned by state cacti with roots in $\mathcal A_u$: Notice that $\mathcal D_{\mathcal A_u}$ is spanned by a disjoint union of cycles and is comprised of a subset of the non-top linked SCCs of $\mathcal D(\bar A)$. For each non-top linked SCC $i$ of $\mathcal{D}_{\mathcal{A}_{u}}$, let $C^{i}$ be a cycle in the $i$-th SCC which contains a state vertex, say $\acute{x}^{i}$, from $\mathcal{A}_{u}$. (Note that such $C^{i}$ exists for all $i$ since $\mathcal{A}_{u}$ contains state vertices from all the non-top linked SCCs.) Hence, another application of Lemma 2 yields that there exist a disjoint collection of chains in $\mathcal{D}_{\mathcal{A}_{u}}$ (as many as the number of non-top linked SCCs in $\mathcal{D}_{\mathcal{A}_{u}}$) such that the first element of each chain, say chain $i$, is the cycle $C^{i}$ and the collection spans $\mathcal{D}_{\mathcal{A}_{u}}$. Now observe that each chain in turn is spanned by a state cacti; specifically, by removing an edge appropriately from the first element (cycle) $C^{i}$ of chain $i$, we may obtain a state cactus in which the state vertices in $C^{i}$ form a state stem with $\acute{x}^{i}\in\mathcal{A}_{u}$ being its root and the remaining cycles in $C^{i}$ are connected from the state stem. Hence, it follows immediately that $\mathcal{D}_{\mathcal{A}_{u}}$ is spanned by a disjoint collection of state cacti with roots in $\mathcal{A}_{u}$.

Hence, there exists a disjoint union of state cacti that spans $\mathcal D(\bar A)$ with roots in $\mathcal{S}_{u}$, which implies by Theorem~1 that $\mathcal S_u$ is a dedicated feasible input configuration. \hfill $\blacksquare$

\textit{Proof of Theorem \ref{Theorem4}} :
From Theorem \ref{FDIC}, it follows that a  (minimal) feasible dedicated input configuration $\mathcal S_u$ must contain a subset $\mathcal U_R$ of right-unmatched vertices w.r.t. some maximum matching of $\mathcal B$ and a subset $\mathcal A_u$ comprising exactly one state variable from each non-top linked SCC of $\mathcal D(\bar A)$. Therefore, for any such pair of subsets $\mathcal{U}_{R}$ and $\mathcal{A}_{u}$, the size of a minimal feasible dedicated input configuration $\mathcal{S}_{u}$ is upper bounded by $|\mathcal{S}_{u}|\leq|\mathcal{U}_{R}|+|\mathcal{A}_{u}|-|\mathcal{U}_{R}\cap\mathcal{A}_{u}|=m+\beta-|\mathcal{U}_{R}\cap\mathcal{A}_{u}|$. The minimum is achieved by maximizing the intersection $|\mathcal{U}_{R}\cap\mathcal{A}_{u}|$; by definition, the maximum value of $|\mathcal{U}_{R}\cap\mathcal{A}_{u}|$ attainable under the given constraints on the sets $\mathcal{U}_{R}$ and $\mathcal{A}_{u}$ is $\alpha$, the maximum top assignability index (see Definition 5). The desired assertion now follows immediately.\hfill $\blacksquare$

\textit{Proof of Theorem \ref{MFDIC}} :
The characterization obtained in Theorem \ref{FDIC} and the development in the proof of Theorem~\ref{Theorem4} implies that there is a one-to-one correspondence between minimal feasible dedicated input configurations and subset pairs $(\mathcal{U}_{R},\mathcal{A}_{u})$ which maximize the intersection $|\mathcal{U}_{R}\cap\mathcal{A}_{u}|$ under the constraints that $\mathcal{U}_{R}$ is  the set of right-unmatched vertices of a maximum matching of $\mathcal{B}$ and $\mathcal{A}_{u}$ is a subset of state vertices with exactly one vertex in each non-top linked SCC of $\mathcal{D}(\bar{A})$. It then readily follows that such maximizing pairs are exactly the ones where $\mathcal{U}_{R}$ corresponds to the set of right-unmatched vertices of a maximum matching of $\mathcal{B}$ with maximum top assignability and $\mathcal{A}_{u}=\mathcal{A}_{R}\cup\mathcal{A}_{u}^{c}$, where $\mathcal{A}_{R}\subset\mathcal{U}_{R}$ denotes the set of state vertices in $\mathcal{U}_{R}$ that belong to non-top linked SCCs, and, $\mathcal{A}_{u}^{c}$ is any subset of state variables comprising only one state
variable from each non-top linked SCC of $\mathcal{D}(\bar{A})$ not assigned by the right-unmatched vertices $\mathcal{U}_{R}$. The desired equivalence follows immediately. \hfill $\blacksquare$

\textit{Proof of Theorem \ref{solP1d}} :
Follows from Theorem \ref{MFDIC}  and by considering duality in linear systems, in other words, by taking $\mathcal D(A^T)$ and designing the corresponding $\bar B$ (using Theorem \ref{MFDIC}), which corresponds to $\bar C^T$.
\hfill $\blacksquare$

\textit{Proof of Theorem \ref{structuralControl}} :
The proof is very similar to that of Theorem 3. We sketch the outline, details are omitted.

[$\Rightarrow$] If  $(\bar A,\bar B)$ is structurally controllable, then by Theorem 1 it follows that it is spanned by a disjoint union of input cacti. By definition, an input cactus is composed by an input stem with edges going from its vertices (either the input or the state vertices) to cycles and/or chains. In addition, an input stem consists of an input with an edge to the root of a state stem, thus, it follows (see proof of Lemma \ref{Lemma2}) that $\mathcal{D}(\bar{A})$ admits a decomposition into state stems and cycles where the root of each state stem is connected from a (distinct) input and the roots of the state stems correspond to a set of right-unmatched vertices $\mathcal{U}^{\prime}_{R}$ associated with a matching of $\mathcal{B}$ (not necessarily maximum). Therefore, by reasoning similarly to Theorem 3, it follows that there exists a subset of right-unmatched vertices $\mathcal U_R$ associated to a maximum matching of $\mathcal B$ that is contained in $\mathcal U_R'$, hence condition i) must hold. Finally, note that each state variable must be reachable from an input (an immediate consequence of the system being spanned by a disjoint union of input cacti) and hence, it follows that in each non-top linked SCC there exists at one state variable that is directly connected from an input.  Thus,  condition ii) must hold.

[$\Leftarrow$] Follows similar arguments as in the proof of Theorem 3 with the following additional notes: (1) by adding an edge  from an input assigned to $\mathcal U_R$ to a state variable in $\mathcal A_u$  that  belongs to the first element (cycle) of a chain, we obtain an input cactus (by the recursive definition of input cactus); and (2)  by Theorem 3, there exists a disjoint union of state cacti with roots in  $\mathcal U_R$ and $\mathcal A_u$ that span $\mathcal{D}(\bar{A})$.
\hfill $\blacksquare$

\textit{Proof of Theorem \ref{sol2a}} :
Follows directly from Theorem \ref{structuralControl} and Theorem \ref{MFDIC} as explained in the main text.\hfill $\blacksquare$

\textit{Proof of Corollary \ref{minsol2a}} :
The corollary follows  from Theorem \ref{sol2a} as explained in the main text.\hfill $\blacksquare$

\textit{Proof of Theorem \ref{solP1}} :
Theorem~\ref{solP1} is an immediate consequence of Theorem \ref{sol2a} and its dual applied to structural output design.\hfill $\blacksquare$

\textit{Proof of Corollary \ref{minsolP1}} :
The corollary is an immediate consequence of  Corollary \ref{minsol2a} and its dual applied to structural output design.\hfill $\blacksquare$

\textit{Proof of Lemma \ref{Lemma5}} :
We start by introducing the notion of augmenting paths: an augmenting path for a matching $M$ is an alternating path (i.e., a sequence of edges alternating between two disjoint sets) with an odd number of edges $e_1e_2 \ldots e_m$ such that $e_{\text{odd}}\in M$ and $e_{\text{even}}\notin M$. In addition, Berge's theorem \cite{Cormen:2001:IA:580470} states that a matching is maximum if and only if it does not have augmenting paths. Now, consider two maximum matchings $M^1,M^2$ of $ \mathcal B\equiv \mathcal B(\mathcal X,\mathcal X,\mathcal E_{\mathcal X,\mathcal X})$ and let $M^1\Delta M^2=M^1\backslash M^2 \cup M^2\backslash M^1$ denote the symmetric difference between the two maximum matchings. In addition, let $\mathcal V^L(M^1\Delta M^2),\mathcal V^D(M^1\Delta M^2)$ denote the left-end vertices and right-end vertices of the edges in $M^1\Delta M^2$ respectively. Take $\mathcal B^{\Delta}\equiv\mathcal B(\mathcal V^L(M^1\Delta M^2),$ $\mathcal V^D(M^1\Delta M^2),M^1\Delta M^2)$, then, by Berge's theorem it follows that no augmenting path (w.r.t. $M^1_{\Delta}=M^1\cap(M^1\Delta M^2)$ or  $M^2_{\Delta}=M^2\cap(M^1\Delta M^2)$ - both maximum matchings of $\mathcal B^{\Delta}$) exists, hence every alternating path w.r.t $M^1_{\Delta}$ (or $M^2_{\Delta}$) in $\mathcal B^{\Delta}$ has even number of edges alternating between $M^1$ and $M^2$. Suppose that $\mathcal V^L(M^1\Delta M^2)$, $\mathcal V^D(M^1\Delta M^2)$ contain $p_1,p_2$ $(\le |\mathcal V^1_R|=|\mathcal V^2_R|=|\mathcal V^1_L|=|\mathcal V^2_L|)$ right and left-unmatched vertices respectively, then there exists $p=p_1+p_2$ alternating paths that start and end in  right/left-unmatched vertices. Let $\mathcal P^L$ denote the collection of $p_2$ alternating paths that start and end in left-unmatched vertices. It follows that $M^*=M^1\Delta \mathcal P^L=M^1\backslash \mathcal P^L\cup \mathcal P^L\backslash M^1$ is a maximum matching of $\mathcal{B}$ (the number of edges is kept the same by the symmetric difference) with unmatched vertices given by $(\mathcal V_R^1,\mathcal V_L^2)$, where $\mathcal V_R^1$ is due to the edges previously in $M^1$ and $\mathcal V_L^2$ induced by $\mathcal P^L$ together with the left-unmatched vertices not in $\mathcal V^L(M^1\Delta M^2)$. \hfill $\blacksquare$

\textit{Proof of Proposition \ref{prop:aux2}} :
\noindent $(1)$ Follows by the construction presented in the main text.

\noindent $(2)$ We provide a constructive proof. Remark that by construction, $\bar{K}\in\mathcal K(M^{u,y},\mathcal S_1,\cdots,\mathcal S_{|M^{u,y}|})$ (see \eqref{defKcal}) yields condition b) in Theorem \ref{Theorem2}. Additionally, if $\mathcal S_i=M^{u,y}$ (for some $i$ and $\mathcal S_i \cap\mathcal S_j=\emptyset$ for all $j\neq i$) then the edge set associated with $\bar K$ is given by $\mathcal F(M^{u,y},\mathcal S_i)$ (defined in \eqref{defFcal}) which leads to a cycle comprising the input-ouput stems described by Proposition \ref{IOdecomp}, therefore it passes through all non-top linked and non-bottom linked SCCs. Hence $\mathcal D(\bar A,\bar B,\bar K,\bar C)$ consists of a single SCC (recall Proposition \ref{prop:aux} which implies that every state vertex can be reached from an effective input and reaches an effective output), which, in addition, contains at least one feedback link (since the edge set of $\bar K$ is non-empty). Thus, condition a) in Theorem \ref{Theorem2} also holds.

\noindent $(3)$ Recall that by construction, the set  $\mathcal K(M^{u,y},\mathcal S_1,\cdots,\mathcal S_{|M^{u,y}|})$,   includes all possible information patterns with exactly $|M^{u,y}|$ feedback links which ensure condition b) in Theorem 2 and, additionally, information patterns with strictly fewer feedback links cannot achieve the design objectives (see the lower bound (12) on the minimum number of feedback edges required and the explanation in the main text). The desired assertion follows immediately.
\hfill $\blacksquare$

\textit{Proof of Theorem \ref{Theorem8}} :
\noindent [$\Leftarrow$] If $\bar B,\bar C$ are constructed as in Corollary~\ref{minsolP1} then the system is structurally controllable/observable and minimality in $\|\bar B\|_0,\|\bar C\|_0$ is achieved. Furthermore,  at least $\max(m,1)$ feedback links, where $m$ denotes the number of right/left-unmatched vertices need to be considered, as obtained in \eqref{lowerboundK}. In addition, the lower bound in \eqref{lowerboundK} can be achieved as described in Proposition \ref{prop:aux2}.

\noindent [$\Rightarrow$] Given the previous developments, it is routine to verify that if any of the conditions on $\bar{B}$, $\bar{C}$, or $\bar{K}$ as described in the hypothesis is not satisfied, the triple cannot be a solution to $\mathcal{P}_{2}$.
\hfill $\blacksquare$

{

\textit{Proof of Theorem \ref{TheoremCorrectness}} :
Step 1 can be performed by executing  depth-first search twice (as explained in chapter 22.5 in \cite{Cormen:2001:IA:580470}, see also Section II), to obtain the DAG and consequently identify the non-top linked SCCs. In particular, the correctness of obtaining the non-top linked SCCs  follows immediately. Further, it incurs in $\mathcal O(|\mathcal X|+|\mathcal E_{\mathcal X,\mathcal X}|)$ complexity. The construction of the weighted bipartite graph $\mathcal B(\mathcal X,\mathcal X\cup\mathcal {I},\mathcal E_{\mathcal X,\mathcal X}\cup \mathcal E_{\mathcal I,\mathcal X})$ follows readily $\mathcal D(\bar A)$ and its DAG representation, which incurs in at most linear complexity. Subsequently,  a minimum weight maximum matching can be obtained using the Hungarian algorithm, which incurs in $\mathcal O(|\mathcal X|^3)$   (see \cite{HungarianAlg}). To see that the set of right-unmatched vertices presented in Step 3 is associated with a  maximum maximum with maximum top-assignability, we notice the following: Let $M'$ be a maximum matching found using the Hungarian algorithm (in particular, it incurs in the minimum cost). Then, $M^*=M'\setminus \{(i,x_j)\in M': i\in \mathcal I \text{ and } x_j \in \mathcal X \}$ is a maximum matching of $\mathcal B\equiv \mathcal B(\mathcal X,\mathcal X,\mathcal E_{\mathcal X,\mathcal X})$, because edges from $\mathcal E_{\mathcal X,\mathcal X}$ incur lesser cost than those in $\mathcal E_{\mathcal I,\mathcal X}$, and the latter are  only used if no edge from $\mathcal E_{\mathcal X,\mathcal X}$ can be used to increase the matching. Additionally, notice that all edges from $\mathcal E_{\mathcal I,\mathcal X}\cap M'$ have one of their endpoints in vertices from $\mathcal X$, which we represent by $\Theta$. Hence, those same vertices become right-unmatched vertices associated with $M^*$. Futhermore, each vertex in $\Theta$ belongs to a different non-top linked SCC, by construction of $\mathcal E_{\mathcal I,\mathcal X}$. In fact, we notice that $\Theta$ comprises the maximum number of state vertices from the set of right-unmatched vertices  in distinct non-top linked SCCs. In other words, the maximum matching $M^*$ has a set of right-unmatched vertices which ensures maximum top-assignability. This last claim is a consequence of noticing that to obtain more right-unmatched vertices in distinct non-top linked SCCs either there exists another another maximum matching $M''$ with an additional edge from $\mathcal E_{\mathcal I,\mathcal X}$ which increases the total cost and consequently $M''\setminus \{(i,x_j)\in M'': i\in \mathcal I \text{ and } x_j \in \mathcal X \}$ is not a maximum matching of $\mathcal B$.  Finally, Step 4 has  linearly computational complexity  and it follows, by Theorem~\ref{MFDIC}, that $\mathcal S_u$ is a minimum feasible dedicated input configuration. \hfill $\blacksquare$
}

\textit{Proof of Corollary \ref{Corollary5}} :
By Theorem \ref{solP1d} and Theorem \ref{solP1} it follows that  solutions to $\mathcal P_1^d$ and $\mathcal P_1$ may be computed by performing simple constructions (without incurring additional complexity) once  a pair of minimal feasible dedicated input/output configurations are available (the latter may  be determined in polynomial complexity in the number of state variables, as stated in Theorem \ref{TheoremCorrectness}). A solution to $\mathcal P_2$ can also be determined by incurring polynomial complexity since it requires a minimal solution to problem $\mathcal P_1$ (may be determined with polynomial complexity) and  a mix-pairing $\bar{K}$ of $\mathcal D(\bar A,\bar B,\bar C)$ (see Theorem~\ref{Theorem8}), the latter can be constructed using a polynomial complexity procedure (also  without increasing the overall complexity), see Remark \ref{rem:mix_pairing}.\hfill $\blacksquare$

\vspace{-0.5cm}

\small
\bibliographystyle{IEEEtran}
\bibliography{IEEEabrv,tac2013}

\end{document}